\documentclass[11pt]{amsart}
\usepackage[english]{babel}
\usepackage{amscd,amssymb,amsmath}
\usepackage{graphicx}
\begin{document}

\title [OPTIMAL QUADRATURE FORMULAS]
{OPTIMAL QUADRATURE FORMULAS IN THE SENSE\\ OF SARD IN
$W_2^{(m,m-1)}$ SPACE}

\author {Kh.M.Shadimetov, A.R.Hayotov}

\address{Kh. M. Shadimetov, A. R. Hayotov\\ Institute of mathematics and information technologies,
Tashkent, Uzbekistan.}
\email {hayotov@mail.ru}

\begin{abstract}
In this paper in $W_2^{(m,m-1)}(0,1)$ space the problem of
construction of optimal quadrature formula in the sense of Sard is
considered and using S.L.Sobolev's method it is obtained new
optimal quadrature formula of such type. For the optimal
coefficients explicit formulas are obtained. Furthermore, the
numerical results which confirm the theoretical results of this
work is given.


\end{abstract}

\maketitle

\section{Introduction. Statement of the problem}

We consider the following quadrature formula
    $$
\int\limits_0^1 {\varphi (x)dx \cong } \sum\limits_{\beta  = 0}^N
{C_\beta  } \varphi (x_\beta  )\eqno     (1.1)
$$
with the error functional
    $$
\ell (x) = \varepsilon _{[0,1]} (x) - \sum\limits_{\beta  = 0}^N
{C_\beta  } \delta (x - x_\beta  ),\eqno   (1.2)
$$
where $C_\beta  $ are the coefficients and $x_\beta$ are the nodes
of formula (1.1), $x_\beta   \in [0,1]$,  $\varepsilon _{[0,1]}
(x)$ is the indicator of the interval $[0,1],$ $\delta (x)$ is
Dirac's delta-function, function $\varphi (x)$ belongs to Hilbert
space $W_2^{(m,m-1)}(0,1)$. The norm of functions in this space is
defined by the following equality
    $$
\left\| {\varphi (x)|W_2^{(m,m-1)}(0,1)} \right\| = \left\{
{\int\limits_0^1 {\left( \varphi^{(m)}(x)+\varphi^{(m-1)}(x)
\right)} ^2 dx} \right\}^{1/2}.\eqno   (1.3)
$$

The difference
    $$
(\ell (x),\varphi (x)) = \int\limits_0^1 {\varphi (x)dx -
\sum\limits_{\beta  = 0}^N {C_\beta  \varphi (x_\beta  )} }  =
\int\limits_{ - \infty }^\infty  {\ell (x)\varphi (x)dx}\eqno
(1.4)
$$
is called \emph{the error} of the quadrature formula (1.1). The
error of the formula (1.1) is linear functional in
$W_2^{(m,m-1)*}(0,1)$, where $W_2^{(m,m-1)*}(0,1)$ is the
conjugate space to the space $W_2^{(m,m-1)}(0,1)$.

By Cauchy-Schwartz inequality
$$
|(\ell(x),\varphi(x))|\leq \|\varphi(x)|W_2^{(m,m-1)}(0,1)\|\cdot
\|\ell(x)|W_2^{(m,m-1)*}(0,1)\|
$$
the error (1.4) of formula (1.1) is estimated with the help of the
norm
$$\left\| {\ell (x)|W_2^{(m,m-1)*}(0,1) } \right\| = \mathop {\sup
}\limits_{\left\| {\varphi (x)|W_2^{(m,m-1)}(0,1)} \right\| = 1}
\left| {\left( {\ell(x) ,\varphi(x) } \right)} \right|
$$
of the error functional (1.2). Consequently, estimation of the
error of the quadrature formula (1.1) on functions of the space
$W_2^{(m,m-1)}(0,1)$ is reduced to finding the norm of the error
functional $\ell(x)$ in the conjugate space $W_2^{(m,m-1)*}(0,1)$.

Obviously the norm of the error functional $\ell (x)$ depends on
the coefficients $C_{\beta}$ and the nodes $x_{\beta}$. The
problem of finding the minimum of the norm of the error functional
$\ell (x)$ by coefficients $C_{\beta}$ and by nodes $x_{\beta}$,
is called the \emph{S.M.Nikolskii problem}, and obtained formula
is called \emph{optimal quadrature formula in the sense of
Nikolskii}. This problem was first  considered by S.M.Nikolskii
[15]. The problem  were further investigated by many authors for
various cases (see e.g. [2-5,16,33] and references therein).
Minimization of the norm of the error functional $\ell(x)$ by
coefficients $C_{\beta}$ when the nodes are fixed is called
\emph{Sard's problem}. And the obtained formula is called
\emph{optimal quadrature formula in the sense of Sard}. This
problem was first investigated by A.Sard [17].

There are several methods for the construction of optimal
quadrature formulas in the sense of Sard such as  the spline
method, the $\varphi(x)-$ function method (see e.g. [2], [24]) and
Sobolev's method. Note the Sobolev method is based on the
construction of a discrete analogue of a linear differential
operator (see e.g. [27-29]). In different spaces based on these
methods, the Sard problem was investigated by many authors (see,
for example, [1,2,4,6-8,11-14,18-20,23-29,31,32] and references
therein).

Furthermore, explicit formulas for coefficients of optimal
quadrature formulas for any $m$ and for any number $N+1$ of the
nodes $x_{\beta}$ in the space $L_2^{(m)}$ when the nodes are
equally spaced were obtained in the works [11,18,19]. Here
$L_2^{(m)}$ is the Sobolev space of functions, which $m-$th
generalized derivative is square integrable. In the work [20]
positiveness of the coefficients of optimal quadrature formulas in
$L_2^{(m)}(0,1)$ space is investigated.

It should be noted that in the works [8,12], using the $\varphi-$
function method, investigated the problem of construction of
optimal quadrature formulas in the sense of Sard which are exact
for solutions of linear differential equations and several
examples given for some number of the nodes.

The main aim of the present paper is to solve the Sard problem in
the space $W_2^{(m,m-1)}(0,1)$ using S.L.Sobolev's method for any number $N+1$
of the nodes $x_{\beta}$,  i.e. finding the coefficients $C_\beta$
satisfying the following equality
    $$
\left\| {\mathop \ell \limits^ \circ  (x)|W_2^{(m,m-1)*}(0,1)} \right\| =
\mathop {\inf }\limits_{C_{_\beta  } } \left\| {\mathop \ell
\limits^{} (x)|W_2^{(m,m-1)*}(0,1)} \right\|. \eqno  (1.5)
$$

Thus, in order to construct the optimal quadrature formula in the
sense of Sard in the space $W_2^{(m,m-1)}(0,1)$ we need
consequently to solve the following problems.

\textbf{Problem 1.} \emph{Find the norm of the error functional
 $\ell(x)$ of quadrature formulas (1.1) in the space $W_2^{(m,m-1)*}(0,1)$.}

\textbf{Problem 2.} \emph{Find the coefficients
$C_\beta $ which satisfy equality (1.5) when the nodes
$x_\beta$ are fixed.}

The structure of the present paper is: in the second section the
extremal function which corresponds to the error functional
$\ell(x)$ is found and with its help  representation of the norm of
the error functional (1.2) is calculated, i.e. the problem 1 is solved;
in the third section in order to find the minimum of the
quantity $\left\| \ell \right\|^2 $ by coefficients
$C_\beta$ the system of linear equations is obtained for the
coefficients of optimal quadrature formula in the sense of Sard in
the space $W_2^{(m,m-1)}(0,1)$, moreover existence and uniqueness of the
solution of this system are proved; in the forth section
explicit formulas for the coefficients of optimal quadrature
formula of the form (1.1) are found, i.e. the problem 2 is solved;
finally, in the fifth section the results of numerical experiments are given.

\section{The extremal function and the representation of the norm\\ of the error functional
$\ell(x)$}

In order to solve the problem 1, i.e. for calculation the norm of
the error functional (1.2) in the space $W_2^{(m,m-1)*}(0,1)$, it
is used the concept of extremal function of given functional. The
function $\psi_{\ell}(x)$ is called \emph{the extremal function}
for the functional $\ell (x)$ (see, [28]), if the following
equality is fulfilled
    $$
\left( {\ell (x),\psi _\ell  (x)} \right) = \left\| {\ell
(x)|W_2^{(m,m-1)*}(0,1) } \right\| \cdot
\left\| {\psi _\ell  (x)|W_2^{(m,m-1)}(0,1)} \right\|.\eqno     (2.1)
$$
This means that the extremal function $\psi_{\ell}(x)$ gives the greatest
error to the difference (1.4) in the space $W_2^{(m,m-1)}(0,1)$.

Since the space $W_2^{(m,m-1)}(0,1)$ is Hilbert space then the extremal
function $\psi_\ell(x)$ in this space is found with the help of
Riesz theorem about general form of a linear continuous functional
on Hilbert spaces. Then for functional $\ell (x)$ and for any
$\varphi(x)\in W_2^{(m,m-1)}(0,1)$ there exists the function $\psi_\ell
(x) \in W_2^{(m,m-1)}(0,1)$ for which the following equation is take placed
$$\left( {\ell (x),\varphi (x)}
\right) = \left\langle {\psi _\ell  (x),\varphi (x)} \right\rangle
, \eqno  (2.2)$$ where
$$\left\langle {\psi _\ell  (x),\varphi
(x)} \right\rangle  = \int\limits_0^1 {\left( {\psi _\ell^{(m)}(x) +
\psi _\ell^{(m-1)}(x)} \right)\left( {\varphi^{(m)}(x) + \varphi^{(m-1)}(x)}
\right)dx}\eqno (2.3)
$$
is the inner product defined in the space $W_2^{(m,m-1)}(0,1)$.

Further, we will solve equation (2.2).

Suppose $\varphi(x)$ belongs to the space $\stackrel{\circ}{C}^{(\infty )}\!\!\!(0,1)$,
where $\stackrel{\circ}{C}^{(\infty)}\!\!\!(0,1)$ is the space of
functions, which are infinity differentiable and finite in the
interval $(0,1)$. Then from (2.3), integrating by parts, we obtain
$$
\left\langle {\psi _\ell  (x),\varphi (x)} \right\rangle  =
(-1)^m\int\limits_0^1 {\left( {\psi _\ell ^{(2m)} (x) - \psi _\ell
^{(2m-2)}(x)}\right)\varphi (x)dx}.\eqno     (2.4)
$$
Keeping in mind (2.4) from (2.2) we get
$$
\psi _\ell ^{(2m)} (x)-\psi _\ell^{(2m-2)}(x)=
(-1)^m\ell (x).\eqno   (2.5)
$$
So, when $\varphi(x) \in \stackrel{\circ}{C} ^{(\infty )}\!\!\!(0,1)$
the extremal function $\psi _\ell  (x)$ is a solution of
equation (2.5). But, we have to find the solution of equation
(2.2) when $\varphi (x) \in W_2^{(m,m-1)}(0,1)$. Since the space
$\stackrel{\circ}{C}^{(\infty )}\!\!\!(0,1)$ is densely in the space
$W_2^{(m,m-1)}(0,1)$, then we can approximate arbitrarily exact functions
of the space $W_2^{(m,m-1)}(0,1)$ by a sequence of functions of the space
$\stackrel{\circ}{C}^{(\infty )}\!\!\!(0,1)$. Next for any $\varphi(x)\in
W_2^{(m,m-1)}(0,1)$ we consider the inner product $\left\langle {\psi
_\ell (x),\varphi (x)} \right\rangle $ and, integrating by parts,
we have
$$
\left\langle {\psi _\ell
(x),\varphi (x)} \right\rangle  = \sum\limits_{s=1}^{m-1}
(-1)^s\varphi^{(m-1-s)}(x)\left({\psi _\ell ^{(m+s)}(x)-\psi
_\ell^{(m+s-2)}(x)} \right)|_{x=0}^{x=1}+
$$
$$
+\varphi^{(m-1)}(x)\left(\psi_{\ell}^{(m)}(x)+\psi_{\ell}^{(m-1)}(x)\right)|_{x=0}^{x=1}
+(-1)^m\int\limits_0^1 {\varphi(x)\left( {\psi_\ell
^{(2m)}(x)-\psi_\ell ^{(2m-2)}(x)}\right)}dx.$$ Hence from
arbitrariness $\varphi(x)$ and uniqueness of the function
$\psi_{\ell}(x)$ (up to the function $e^{-x}$ and polynomial of
degree $m-2$), taking into account (2.5), it must be fulfilled the
following equations
$$
\psi _\ell ^{(2m)}(x) -\psi _\ell^{(2m-2)}(x)=(-1)^m\ell (x), \eqno  (2.6)
$$
$$\left( {\psi _\ell^{(m+s)}(x)-\psi _\ell^{(m+s-2)}(x)}\right)|_{x = 0}^{x = 1}  =
 0,\eqno
 (2.7)
 $$
$$\left( {\psi _\ell^{(m)}(x) + \psi _\ell^{(m-1)}(x)}\right)|_{x = 0}^{x = 1}  = 0.
\eqno (2.8)
   $$
Thus, we conclude that the extremal function $\psi_{\ell}(x)$ is
the solution of the boundary value problem (2.6)-(2.8).

The following holds

\textbf{Theorem 2.1.} {\it The solution of the boundary value
problem (2.6)-(2.8) is the extremal function $\psi_{\ell}(x)$ of
the error functional $\ell (x)$ and has the following form
$$
\psi_\ell(x)=(-1)^m\ell (x)*G(x)+ P_{m-2}(x) + d e^{- x},
$$
where
    $$
G(x) = {{{\rm{sign}}x} \over 2}\left(\frac{e^x-e^{-x}}{2}-\sum\limits_{k=1}^{m-1}
\frac{x^{2k-1}}{(2k-1)!}\right) \eqno     (2.9)
$$
is a solution of the equation
  $$
\psi ^{(2m)} (x)-\psi^{(2m-2)}(x)=
\delta (x)\eqno   (2.10)
$$
$d$ is any real number, $P_{m-2}(x)$ is a polynomial of degree $m-2$.}

\textbf{Proof.} It is known that general solution of nonhomogeneous
differential equation consists on sum of a partial solution of
nonhomogeneous  differential equation and general solution of
corresponding homogeneous differential equation.

The homogeneous equation for differential equation (2.6) has the
form
    $$
\psi _\ell ^{(2m)} (x) - \psi _\ell^{(2m-2)}(x) =
0.\eqno (2.11)
$$
It is easy to show that general solution of
equation (2.11) is
$$
P_{2m-3}(x)+d_1e^x+d_2e^{-x}.\eqno   (2.12)
$$
It is not difficult to verify that a partial solution of the
differential equation (2.6) is $$(-1)^m\ell (x)*G(x),$$ where $G(x)$ is
a fundamental solution of equation (2.6) and is defined by
(2.9) and is a solution of equation (2.10), * is
operation of convolution, where convolution of two functions is defined as
$$
f(x)*g(x) = \int\limits_{ - \infty }^\infty  {f(x - y)g(y)dy = }
\int\limits_{ - \infty }^\infty  {f(y)g(x - y)dy}.
$$

The rule of finding a fundamental solution of a linear differential
operator
    $$
L \equiv {{d^n } \over {dx^n }} + a_1 {{d^{n - 1} } \over {dx^{n -
1} }} + ... + a_n ,
$$
where $a_j$ are constants, is given in [30, p.88]. Using this
rule, it is found the function $G(x)$, which is a fundamental
solution of the operator ${{d^{2m}} \over {dx^{2m}}} - {{d^{2m-2}}
\over {dx^{2m-2}}}$ and have the form (2.9).

Thus, we have the following general solution of the equation (2.6)
$$
\psi _\ell(x)=(-1)^m\ell(x)*G(x) + P_{2m-3}(x)+d_1e^x+d_2^{-x},  \eqno     (2.13)
$$
where $P_{2m-3}(x)=b_{2m-3}x^{2m-3}+...+b_1x+b_0$ is a polynomial of degree $2m-3$,
$d_1,\ d_2$ are constants.

In order that in the space $W_2^{(m,m-1)}(0,1)$ the function $\psi
_\ell (x)$ will be unique  (up to the function $e^{-x}$ and
polynomial of degree $m-2$), it has to satisfies the conditions
(2.7), (2.8). Here derivative is in generalized sense and
$$
\psi _\ell ^{(k)}(x)=(-1)^m \ell(x)*G^{(k)}(x) + P_{2m - 3}^{(k)} (x) + d_1 e^x  + ( - 1)^k d_2 e^{ - x} ,
$$
$$
k = 1,2,...,2m - 1,
$$
where
$$
G^{(k)}(x)={{{\rm{sign}}x} \over 2}\left\{
\begin{array}{l}
\displaystyle{{e^x  - e^{ - x} } \over 2} - {{x^{2m - 3 - k} } \over {(2m - 3 - k)!}} - ... - {{x^3 } \over {3!}} - x \mbox{ when } k\mbox{ is even},\\
\displaystyle{{e^x  + e^{ - x} }
\over 2} - {{x^{2m - 3 - k} } \over {(2m - 3 - k)!}} - ... - {{x^2 } \over {2!}} - 1\mbox{ when } k\mbox{ is odd.}\\
\end{array}
\right.
$$
From conditions (2.7) for $s = m - 1$, taking into account
$$
\psi_\ell^{(2m - 1)}(x)=( - 1)^m \ell(x)*\left\{ {{\rm{sign}}x
\cdot{{e^x  + e^{ -x} } \over 4}} \right\} + d_1 e^x  - d_1 e^{ - x} ,
$$
$$
\psi _\ell ^{(2m - 3)} (x) = ( - 1)^m \ell(x)*\left\{ {{{{\rm{sign}}x} \over 2}\left[ {{{e^x  + e^{ - x} } \over 2} - 1} \right]} \right\} + (2m - 3)!b_{2m - 3}  + d_1 e^x  - d_1 e^{ - x} ,
$$
we get
$$
\psi_\ell ^{(2m - 1)}(x)-\psi _\ell ^{(2m - 3)} (x) =
( - 1)^m \ell(x)*{{{\rm{sign}}x} \over 2} - (2m - 3)!b_{2m - 3}  =
$$
$$
=( - 1)^m \left( {\ell(y),{{{\rm{sign}}(x - y)} \over 2}} \right) - (2m - 3)!b_{2m - 3} .
$$
Hence for $x = 0$:
$$
\psi _\ell ^{(2m - 1)}(0) - \psi _\ell ^{(2m - 3)}(0)=
( - 1)^m \left( {\ell(y),{{{\rm{sign}}( - y)} \over 2}} \right) - (2m - 3)!b_{2m - 3}  =
$$
$$
=-{{( - 1)^m } \over 2}\left( {\ell (y),1} \right)-(2m - 3)!b_{2m - 3}  = 0,
$$
for $x = 1$:
$$
\psi _\ell ^{(2m - 1)} (1) - \psi _\ell ^{(2m - 3)}(1)=( - 1)^m
\left( {\ell(y),{{{\rm{sign}}(1 - y)} \over 2}} \right) - (2m - 3)!b_{2m - 3}  =
$$
$$
= {{( - 1)^m } \over 2}\left( {\ell(y),1} \right) - (2m - 3)!b_{2m - 3}  = 0.
$$
Then
$$
b_{2m - 3}=0\,\,\,\,\,\mbox{ and }\,\,\,\,(\ell(y),1) = 0.\eqno     (2.14)
$$
When $s = m - 2$  from conditions (2.7), taking into account
$$
\psi _\ell ^{(2m - 2)} (x)=( - 1)^m \ell(x)*
\left\{ {{\rm{sign}}x\cdot {{e^x  - e^{ - x} } \over 4}} \right\} + d_1 e^x  + d_1 e^{ - x} ,
$$
$$
\psi _\ell ^{(2m - 4)}(x)=( - 1)^m \ell _N (x)*
\left\{ {{{{\rm{sign}}x} \over 2}\left[ {{{e^x  - e^{ - x} } \over 2} - x} \right]} \right\} +
$$
$$
+(2m - 3)!b_{2m - 3} x + (2m - 4)!b_{2m - 4}  + d_1 e^x  + d_2 e^{ - x},
$$
we have
$$
\psi_\ell^{(2m - 2)}(x)-\psi _\ell^{(2m - 4)}(x)
=( - 1)^m \ell(x)*{{x \cdot {\rm{sign}}x} \over 2}-
(2m-3)!b_{2m - 3} x - (2m - 4)!b_{2m - 4}.
$$
Hence, taking account of (2.14), for $x=0$:
$$
(-1)^m\left({\ell(y),{{0 - y}\over 2}{\rm{sign}}(0 - y)}\right)-(2m - 4)!b_{2m - 4}=
$$
$$
={{( - 1)^m } \over 2}\left( {\ell(y),y} \right) - (2m - 4)!b_{2m - 4}  = 0,
$$
for $x = 1$:
$$
(-1)^m \left({\ell(y),{{1 - y}\over 2}{\rm{sign}}(1 - y)}\right)-(2m - 4)!b_{2m - 4}  =
$$
$$
=-{{(-1)^m }\over 2}\left( {\ell(y),y} \right) - (2m - 4)!b_{2m - 4}  = 0.
$$
Hence
$$
b_{2m - 4}  = 0\mbox{ and }(\ell(y),y) = 0,\eqno (2.15)
$$
and so on, continuing by this manner,  for $s=m-3,m-4,...,2,1$ we obtain
$$
b_{m-2+s}=0,\ \ \ (\ell(y),y^{m-1-s})=0.\eqno (2.16)
$$
Combining (2.14), (2.15), (2.16) we get
$$
b_{2m-3-s}=0,\ \ \  s=0,1,...,m-2,\eqno (2.17)
$$
$$
(\ell(y),y^s),\ \ \ s=0,1,...,m-2.\eqno (2.18)
$$
From the condition (2.8), keeping in mind (2.17), we have
$$
\psi_\ell ^{(m)}(x)=( - 1)^m \ell(x)*G^{(m)} (x) +
d_1 e^x  + ( - 1)^m d_2 e^{ - x},
$$
$$
\psi _\ell ^{(m - 1)}(x)=( - 1)^m \ell (x)*G^{(m - 1)}(x) +
d_1 e^x  + ( - 1)^{m - 1} d_2 e^{ - x}.
$$
Hence
$$
\psi _\ell ^{(m)} (x) + \psi _\ell ^{(m - 1)} (x) =
$$
$$
=( - 1)^m \ell (x)*\left\{{{{\mathrm{sign}x} \over 2}
\left( {e^x  - {{x^{m - 2} } \over {(m - 2)!}} -
{{x^{m - 3} } \over {(m - 3)!}} - ... - x - 1} \right)} \right\} + d_1 e^x .
$$
For $x = 0$, taking into account (2.18),
$$
\psi _\ell ^{(m)}(0) + \psi _\ell ^{(m - 1)}(0)=
(-1)^m \left( {\ell(y),{{{\rm{sign}}(0 - y)} \over 2}}
\right. \cdot \left( {e^{ - y}  - {{( - y)^{m - 2} } \over {(m - 2)!}}} \right. -
$$
$$
- \left. {\left. {{{( - y)^{m - 3} } \over {(m - 3)!}} - ...
+ y - 1} \right)} \right) + d_1  =  - {{( - 1)^m } \over 2}
\left( {\ell _N (y),e^{ - y} } \right) + d_1  = 0.
$$
For $x = 1$:
$$
\psi_\ell^{(m)}(1)+\psi_\ell^{(m - 1)}(1)=
(-1)^m \left( {\ell(y),{{{\rm{sign}}(1 - y)}\over 2}}\right.
\cdot \left( {e^{1 - y}  - {{(1 - y)^{m - 2} } \over {(m - 2)!}}} \right. -
$$
$$
-\left.{\left. {{{(1 - y)^{m - 3}}\over {(m - 3)!}} -...-(1 - y)-1}\right)}
\right) + d_1 e = {{( - 1)^m e} \over 2}\left( {\ell _N (y),e^{ - y} } \right) + d_1 e = 0.
$$
Whence
$$
d_1  = 0,\eqno     (2.19)
$$
$$
\left({\ell(y),e^{ - y} } \right) = 0.\eqno     (2.20)
$$
Taking into account equalities (2.17)- (2.20) and denoting $d_2  =
d$, we get the statement of the theorem. Theorem 2.1 is
proved.\hfill $\Box$

The equalities (2.18) and (2.20) mean that our quadrature formula
is exact for the function $e^{-x}$ and for any polynomial of
degree up to $m-2$.

Now, using the result of theorem 2.1, we immediately  obtain the
representation of square of the norm of the error functional (1.2)
$$
\left\| \ell(x) |W_2^{(m,m-1)*}(0,1)\right\|^2 = (\ell(x),\psi _\ell
(x)) = (-1)^m\Bigg[\sum\limits_{\beta  = 0}^N {\sum\limits_{\gamma  = 0}^N
{C_\beta  C_\gamma  \,G(x_\beta   - x_\gamma  ) - } }
$$
 $$
- 2\sum\limits_{\beta  = 0}^N C_\beta  \int\limits_0^1 {G(x -
x_\beta  )dx + \int\limits_0^1 {\int\limits_0^1 {G(x - y)dxdy} } }\Bigg],
\eqno   (2.21)
$$

Thus the problem 1 is solved.

Further in sections 3 and 4 we solve the problem 2.

\section{The system for the coefficients of optimal quadrature formulas in the
space $W_2^{(m,m-1)}(0,1)$}

Assume that the nodes $x_\beta$  of the quadrature formula (1.1)
are fixed. The error functional (1.2) satisfies conditions (2.18)
and (2.20). The norm of the error functional $\ell (x)$ is
multidimensional function with respect to the coefficients
$C_\beta$ $(\beta  = \overline {0,N} )$. For finding the point of
conditional minimum of the expression (2.21) under the conditions
(2.18) and (2.20) we apply the method of undetermined multipliers
of Lagrange.

We denote $\mathbf{C}=(C_0,C_1,...,C_N)$ and
$\lambda=(\lambda_0,\lambda_1,...,\lambda_{m-1})$.

Consider the function
$$
\Psi (\mathbf{C},\lambda) =\left\| \ell(x) \right\|^2  -
2(-1)^m \sum\limits_{\alpha=0}^{m-2}\lambda_{\alpha}
\left({\ell(x),x^{\alpha}}\right)- 2(-1)^m \lambda_{m-1}
\left(\ell(x),e^{-x}\right).
$$
Equating to 0 partial derivatives of $\Psi(\mathbf{C},\lambda)$ by $C_\beta$
$(\beta =\overline{0,N})$
and $\lambda_0,\lambda_1,...,\lambda_{m-1}$, we
get the following linear system
$$
\sum\limits_{\gamma=0}^N C_\gamma G(x_\beta   - x_\gamma  ) +
\sum\limits_{\alpha=0}^{m-2}\lambda_{\alpha}x_{\beta}^{\alpha}+ \lambda_{m-1}e^{-x_{\beta}}
= f_m(x_{\beta}),\,\,\,\,\beta  = \overline {0,N},\eqno     (3.1)
$$
$$
\sum\limits_{\gamma=0}^N C_\gamma x_{\gamma}^{\alpha}=\frac{1}{\alpha+1},\
\alpha=0,1,...,m-2,\eqno     (3.2)
$$
$$
\sum\limits_{\gamma=0}^NC_\gamma e^{-x_{\gamma}}=1-e^{-1},  \eqno (3.3)
$$
where $G(x)$ is defined by equality (2.9),
    $$
f_m(x_\beta ) = \int_0^1 {G(x - x_\beta  )dx}. \eqno  (3.4)
$$

The system (3.1)-(3.3) has a unique solution and this solution
gives minimum to $\left\| \ell(x)\right\|^2 $ under the conditions
(3.2), (3.3). Uniqueness of the solution of such type of systems
is discussed in [28, 29].

It should be noted that existence and uniqueness of optimal
quadrature formulas in the sense of Sard is also investigated in
[12].

Now in (2.21) we will do change of variables $C_\beta = \overline
C _\beta + C_{1\beta }$. Then (2.21) and the system (3.1)-(3.3)
have the following form:
$$
\left\| \ell  \right\|^2  = (-1)^m\Bigg[\sum\limits_{\beta  = 0}^N
{\sum\limits_{\gamma  = 0}^N {\overline C _\beta  \overline C
_\gamma  G(x_\beta   - x_\gamma  )}  - 2\sum\limits_{\beta  = 0}^N
{(\overline C _\beta   + C_{1,\beta } )\int\limits_0^1 {G(x -
x_\beta  )dx + } } }
$$
$$
+ \sum\limits_{\beta  = 0}^N {\sum\limits_{\gamma  = 0}^N {\left(
{2\overline C _\beta  C_{1,\gamma }  + C_{1,\beta } C_{1,\gamma }
} \right)G(x_\beta   - x_\gamma  )}  + \int\limits_0^1
{\int\limits_0^1 {G(x - y)dxdy} } }\Bigg],\eqno   (3.5)
$$
  $$
\sum\limits_{\gamma  = 0}^N {\overline C _\gamma  G(x_\beta   -
x_\gamma  ) + \sum\limits_{\alpha=0}^{m-2}\lambda_{\alpha}x_{\beta}^{\alpha}+
\lambda_{m-1}e^{-x_{\beta}} =
F_m(x_\beta),\,\,\,\,\beta  = \overline {0,N} } ,\eqno     (3.6)
$$
$$
\sum\limits_{\gamma  = 0}^N \overline C _\gamma  x_\gamma^{\alpha}
= 0,\ \ \ \alpha=\overline{0,m-2},\eqno    (3.7)
$$
$$
\sum\limits_{\gamma=0}^N{\overline C _\gamma e^{-x_\gamma}
= 0}, \eqno  (3.8)
$$
where  $F_m(x_\beta  ) = f_m(x_\beta  ) - \sum\limits_{\gamma = 0}^N
{C_{1,\gamma }G(x_\beta   - x_\gamma  )} $, $C_{1\beta } $ is a
partial solution of equations (3.2), (3.3).

Hence we directly get that the minimization of (2.21) under the
conditions (3.2), (3.3) with respect to $C_\beta$ is equivalent to the minimization
of expression (3.5) with respect to $\overline C _\beta$ under the conditions
(3.7), (3.8). Therefore it is sufficient to prove that the
system (3.6)-(3.8) has a unique solution with respect to unknowns
$\overline{\mathbf{C}}  = (\overline C _0 ,\overline C _1
,...,\overline C _N ),$ $\lambda = (\lambda_0,\lambda_1,...,\lambda_{m-1})$ and this
solution gives conditional minimum to the expression $\left\| \ell
\right\|^2 $.

From the theory of conditional extremum it is known the sufficient
condition in which the solution of the system (3.6) - (3.8) gives
conditional minimum to the expression $\left\| \ell \right\|^2$ on
the manifold (3.7), (3.8). It consists on positiveness of the
quadratic form
$$
\Phi (\overline{\mathbf{C}} ) = \sum\limits_{\beta  = 0}^N
{\sum\limits_{\gamma  = 0}^N {{{\partial ^2 \Psi } \over {\partial
\overline C _\beta  \partial \overline C _\gamma  }}} } \overline
C _\beta  \overline C _\gamma  \eqno   (3.9)
$$
on the set of the vectors  $\overline{\mathbf{C}}  = (\overline C
_0,\overline C _1 ,...,\overline C _N )$ under the condition
$$
S\overline{\mathbf{C}}  = 0,\eqno     (3.10)
$$
where $S$ is the following matrix of equations (3.7), (3.8):
$$
S = \left( \begin{array}{cccc}
   1 & 1  &  \cdots &    1 \\
 x_0 & x_1&  \cdots & x_N  \\
 \vdots &\vdots &\ddots &\vdots \\
x_0^{m-2} & x_1^{m-2}&  \cdots & x_N^{m-2}  \\
e^{-x_0} & e^{-x_1}&  \cdots & e^{-x_N}  \\
   \end{array}
\right).
$$

We will show that in our case this condition is fulfilled.

\textbf{Theorem 3.1.} {\it For any nonzero vector $
\overline{\mathbf{C}} \in R^{N + 1}$ lying in the subspace
$S\overline {\mathbf{C}} = 0$, the function $\Phi (\overline
{\mathbf{C}} )$ is strictly positive.}

{\bf Proof.} Using the definition of the function  $\Psi
(\mathbf{C},\lambda)$ and equations (3.5), (3.7), (3.8) from
(3.9) we get
    $$
\Phi \left( {\overline {\mathbf{C}} } \right) =
2(-1)^m\sum\limits_{\beta = 0}^N {\sum\limits_{\gamma  = 0}^N {G(x_\beta
- x_\gamma )\overline C _\beta  } } \overline C _\gamma  .\eqno
(3.11)
$$

Consider the linear combination of delta functions
    $$
\delta _{\overline{\mathbf{C}}}(x) = \sqrt 2 \sum\limits_{\beta =
0}^N {\overline C _\beta  \delta (x - x_\beta  ).}\eqno     (3.12)
$$

By virtue of the condition (3.10) this functional belongs to the
space $W_2^{(m,m-1)*}(0,1)$. So, it has the extremal function $u_{\overline
{\mathbf{C}}}(x) \in W_2^{(m,m-1)}(0,1)$ which is a solution of the
equation
 $$
\left({d^{2m} \over dx^{2m}}-{d^{2m-2} \over dx^{2m-2}}\right)
u_{\overline{\mathbf{C}}}(x) =(-1)^m \delta_{\overline
{\mathbf{C}}}(x).\eqno (3.13)
$$
As $u_{\overline {\mathbf{C}} } (x)$ we can take linear
combination of shifts of the fundamental solution $G(x)$:
$$
u_{\overline {\mathbf{C}}} (x) = \sqrt 2 (-1)^m \sum\limits_{\beta  =
0}^N {\overline C _\beta  G(x - x_\beta  )}.
$$
The square of its norm in the space $W_2^{(m,m-1)}(0,1)$ coincide with $\Phi
\left( {\overline{\mathbf{C}}} \right)$, i.e.
$$
\left\| {u_{\overline {\mathbf{C}} } (x)|W_2^{(m,m-1)}(0,1)} \right\|^2  =
\left( {\delta _{\overline {\mathbf{C}}} (x),u_{\overline
{\mathbf{C}} } (x)} \right) = 2(-1)^m\sum\limits_{\beta  = 0}^N
{\sum\limits_{\gamma  = 0}^N {\overline C _\beta  \overline C
_\gamma  G(x_\beta  } }  - x_\gamma  ).
$$
Hence clearly that for nonzero $\overline {\mathbf{C}}$ the
function $\Phi \left( {\overline {\mathbf{C}}} \right) $ is
strictly positive.

Theorem 3.1 is proved.\hfill  $\Box$ \\[1mm]

If the nodes $x_0 ,x_1 ,...,x_N$  are selected such that the matrix
$S$ has right inverse matrix, then the system (3.6) - (3.8) has
a unique solution. Then the system (3.1)-(3.3) also has a unique
solution.

\textbf{Theorem 3.2.}  \emph{If the matrix $S$ has right inverse
matrix, then the main matrix $Q$ of the system (3.6) - (3.8) is
nonsingular.}

{\bf Proof.} We denote by $M$ the matrix of quadratic form $\frac{(-1)^{m}}{2}\Phi
\left( {\overline {\mathbf{C}}} \right)$, where $ \Phi \left(
{\overline{\mathbf{C}}} \right)$ is defined by equality (3.11). It is
known that if homogenous system of linear equations has only trivial
solution then corresponding nonhomogeneous system has a unique
solution. Consider homogeneous system corresponding to the system
(3.6) -(3.8) in the following matrix form:
$$
Q\left(
\begin{array}{c}
   \overline{\mathbf{C}}  \\
   \lambda\\
\end{array}
\right) =\left(
\begin{array}{cc}
   M & S^*  \\
   S & 0  \\
\end{array}
 \right)
\left(
\begin{array}{c}
\overline{\mathbf{C}}  \\
\lambda  \\
\end{array}
\right)=0,
\eqno (3.14)
$$
where $S^*$ is the transposed matrix to the matrix $S$.

We verify, that unique solution of (3.14) is identical zero. Suppose
$\overline{\mathbf{C}},\,\,\lambda $ is the solution of
(3.14). Consider the function $\delta_{\overline{\mathbf{C}}}(x)$,
which is determined by equality (3.12). As the extremal function for
the function $\delta _{\overline{\mathbf{C}}}(x)$ we take
the following function:
$$
u_{\overline {\mathbf{C}} } (x) = \sqrt 2(-1)^m\left( \sum\limits_{\beta  =
0}^N {\overline C _\beta  G(x - x_\beta  ) +\sum\limits_{\alpha=0}^{m-2}\lambda_{\alpha}
x^{\alpha}+\lambda_{m-1}e^{-x}}\right).
$$
This is possible because the function $u_{\overline {\mathbf{C}}}(x)$ belongs
to the space $W_2^{(m,m-1)}(0,1)$ and is a solution of equation (3.13).
First $N + 1$ equations of the system (3.14) mean that
$u_{\overline {\mathbf{C}} }(x)$ takes 0 values at all nodes
$x_\beta$. Then for the norm of the functional $\delta
_{\overline{\mathbf{C}}}(x)$ in $W_2^{(m,m-1)*}(0,1)$ we have
$$
\left\|{\delta _{\overline {\mathbf{C}}} (x)|W_2^{(m,m-1)*}(0,1)}
\right\|^2  = \left( {\delta _{\overline {\mathbf{C}} }
(x),u_{\overline {\mathbf{C}} } (x)} \right) = 2(-1)^m
\sum\limits_{\beta  = 0}^N {\overline C _\beta \,u_{\overline
{\mathbf{C}} } (x_\beta  ) = 0}\eqno (3.15)
$$
On the other hand taking into account equations (3.7), (3.8) we get
$$
\left\|\delta _{\overline {\mathbf{C}}}(x)|W_2^{(m,m-1)*}(0,1)
\right\|^2  = \left(\delta_{\overline {\mathbf{C}}}
(x),u_{\overline {\mathbf{C}} } (x) \right) =2(-1)^m\sum\limits_{\beta  = 0}^N
\sum\limits_{\gamma  = 0}^N \overline C _\beta \overline C
_\gamma  G(x_\beta-x_\gamma).\eqno (3.16)
$$
From (3.16) we conclude that (3.15) is possible if $\overline{\mathbf{C}}=0$.
Then from first $N + 1$ equations of the system (3.14) we obtain
$$
S^* \lambda= 0.\eqno     (3.17)
$$
By assertion of the theorem, the matrix $S$ has right inverse
matrix, then $S^*$ has left inverse matrix. Then from (3.17) we
conclude that the solution $\lambda$ also is equal to zero.

Theorem 3.2 is proved.\hfill $\Box$\\[1mm]

From (2.21) and theorems 3.1, 3.2 it follows that in fixed values
of the nodes $x_\beta$ the square of the norm of the error functional
$\ell(x)$ being quadratic functions of the coefficients $C_\beta$
has a unique minimum in some concrete value $C_\beta= \mathop
{C_\beta }\limits^ \circ$.

As said in the first section the quadrature formulas with the coefficients
$\stackrel{\circ}{C}_\beta \,\,(\beta  = \overline {0,N} )$,
corresponding to this minimum in fixed nodes $x_\beta$ is called
\emph{optimal quadrature formula in the sense of Sard} and
$\stackrel{\circ}{C}_\beta \,\,(\beta = \overline {0,N} )$ are
called \emph{optimal coefficients}.

Below for convenience the optimal coefficients
$\stackrel{\circ}{C}_\beta$ remain as $C_\beta$.

\section{Coefficients of optimal quadrature formula in the sense of Sard}

In the present section we will solve the system (3.1)-(3.3) and
will find explicit formula for the coefficients $C_\beta$. Here we
will use similar method suggested by S.L.Sobolev [27] for finding
the coefficients of optimal quadrature formulas in the space
$L_2^{(m)}(0,1)$. Here mainly is used the concept of functions of
discrete argument and operations on them. Theory of discrete
argument functions is given in [28,29]. For completeness we give
some definitions about functions of discrete argument.

Assume that the nodes $x_\beta$ are equal spaced, i.e. $x_\beta=
h\beta,$ $h = {1 \over N}$, $N = 1,2,...$.

Suppose that $\varphi (x)$ and $\psi (x)$ are real-valued functions
of real variable and are defined in real line $\mathbb{R}$.

\textbf{Definition 4.1.} A function $\varphi (h\beta )$ is called
\emph{function of discrete argument} if it is given on some set
of integer values of $\beta$.

\textbf{Definition 4.2.} \emph{By inner product} of two discrete
functions $\varphi(h\beta )$ and $\psi (h\beta )$ is called the
number
    $$
\left[ {\varphi(h\beta),\psi(h\beta) } \right] = \sum\limits_{\beta  =  - \infty
}^\infty  {\varphi (h\beta ) \cdot \psi (h\beta )},
$$
if the series on the right hand side of the last equality converges
absolutely.

\textbf{Definition 4.3.} \textit{By convolution} of two discrete functions
$\varphi(h\beta )$ and $\psi (h\beta )$ is called the inner
product
$$
\varphi (h\beta )*\psi (h\beta ) = \left[ {\varphi (h\gamma
),\psi (h\beta  - h\gamma )} \right] = \sum\limits_{\gamma  =  -
\infty }^\infty  {\varphi (h\gamma ) \cdot \psi (h\beta  - h\gamma
)}.
$$

Now we turn on to our problem.

Suppose that $C_\beta=0$  when $\beta  < 0$ and $\beta  > N$.
Using above mentioned definitions the system (3.1)-(3.3) we
rewrite in the convolution form
    $$
G(h\beta )*C_\beta   + P_{m-2}(h\beta) + de^{-h\beta} =
f_m(h\beta ),\,\,\,\,\beta  = 0,1,...,N\eqno     (4.1)
$$
$$
C_\beta=0,\,\mbox{ when }\,\,\,\,\,\beta  < 0\mbox{ and }\beta
> N,\eqno (4.2)
$$
$$
\sum\limits_{\beta=0}^N C_\beta\cdot
(h\beta)^{\alpha}=\frac{1}{\alpha+1},\ \ \  \alpha=0,1,...,m-2,
\eqno (4.3)
$$
$$
\sum\limits_{\beta=0}^N C_\beta\cdot e^{-h\beta} = 1-e^{-1},\eqno (4.4)
$$
where
$$
f_m(h\beta)=\frac{e^{h\beta}+e^{-h\beta}+e^{1-h\beta}+e^{h\beta-1}-4}{4}-
\sum\limits_{k=1}^{m-1}\frac{(h\beta)^{2k}+(1-h\beta)^{2k}}{2\cdot
(2k)!}.\eqno         (4.5)
$$

Consider the following problem.

\textbf{Problem A.} \emph{Find the discrete function $C_\beta$,
polynomial $P_{m-2}(h\beta)$ of degree $m-2$ and unknown constant
$d$ which satisfy the system (4.1)-(4.4) for given $f_m(h\beta
)$.}

Further we investigate the problem A and instead of $C_\beta$ we introduce the functions
$$
v(h\beta ) = G(h\beta )*C_\beta\eqno     (4.6)
$$
and
$$
u(h\beta)=v\left({h\beta}\right)+P_{m-2}(h\beta) + d\ e^{-h\beta}.\eqno     (4.7)
$$
In such statement it is necessary to express the coefficients $C_\beta$ by the function
$u(h\beta )$. For this we have to construct such operator
$D_m(h\beta)$ which satisfies the equality
    $$
D_m(h\beta )*G(h\beta ) =\delta(h\beta) ,\eqno    (4.8)
$$
where $\delta (h\beta )$ is equal to 0 when $\beta  \ne 0$ and is
equal to 1 when $\beta  = 0$, i.e. $\delta(h\beta)$ is the
discrete delta-function.

In connection with this in [21,22] the discrete analogue
$D_m(h\beta )$ of the  operator ${{d^{2m}} \over {dx^{2m}
}}-{{d^{2m-2}} \over {dx^{2m-2} }}$, which satisfies equation
(4.8) is constructed and its some properties are investigated.

The following theorems are proved in the works [21,22].

\textbf{Theorem 4.1.} {\it The discrete analogue of the
differential operator ${{d^{2m}} \over {dx^{2m} }}-{{d^{2m-2}} \over {dx^{2m-2} }}$
satisfying the equation (4.8) has the form
$$
D_m(h\beta)=\frac{1}{p_{2m-2}^{(2m-2)}}\left\{
\begin{array}{ll}
\sum\limits_{k=1}^{m-1}A_k \lambda_k^{|\beta|-1},& |\beta|\geq 2,\\
-2e^h+\sum\limits_{k=1}^{m-1}A_k,& |\beta|=1,\\
2C+\sum\limits_{k=1}^{m-1}\frac{A_k}{\lambda_k},&
\beta=0,
\end{array}
\right.\eqno (4.9)
$$
where
$$
C=1+(2m-2)e^h+e^{2h}+\frac{e^h\cdot p_{2m-3}^{(2m-2)}}{p_{2m-2}^{(2m-2)}},\eqno (4.10)
$$
$$
A_k=\frac{2(1-\lambda_k)^{2m-2}[\lambda_k(e^{2h}+1)-e^h(\lambda_k^2+1)]p_{2m-2}^{(2m-2)}}
{\lambda_k\ P_{2m-2}'(\lambda_k)},\eqno (4.11)
$$
$$
\mathcal{P}_{2m-2}(\lambda)=\sum\limits_{s=0}^{2m-2}p_s^{(2m-2)}\lambda^s=(1-e^{2h})(1-\lambda)^{2m-2}-
2(\lambda(e^{2h}+1)-e^h(\lambda^2+1))\times
$$
$$
\times \left[h(1-\lambda)^{2m-4}+\frac{h^3(1-\lambda)^{2m-6}}{3!}E_2(\lambda)+...+
\frac{h^{2m-3}E_{2m-4}(\lambda)}{(2m-3)!}\right],\eqno (4.12)
$$
$p_{2m-2}^{(2m-2)},\ p_{2m-3}^{(2m-2)}$ are the coefficients of
the polynomial $\mathcal{P}_{2m-2}(\lambda)$ defined  by equality
(4.12), $\lambda_k$ are the roots of the polynomial
$\mathcal{P}_{2m-2}(\lambda)$ which absolute values less than 1,
$E_k(\lambda)$ is the Euler-Frobenius polynomial of degree $k$
(see [29]).}
\\[0.2cm]

\textbf{Theorem 4.2.} {\it The discrete analogue $D_m(h\beta)$ of
the differential operator ${{d^{2m}} \over {dx^{2m} }}-{{d^{2m-2}} \over {dx^{2m-2} }}$
satisfies the following equalities

1) $D_m(h\beta)*e^{h\beta}=0,$

2) $D_m(h\beta)*e^{-h\beta}=0,$

3) $D_m(h\beta)*(h\beta)^n=0,$ $n\leq 2m-3$,

4) $D_m(h\beta)*G(h\beta)=\delta(h\beta),$\\
here $G(h\beta)$ is the function of discrete argument
corresponding to the function $G(x)$ defined by equality (2.9)
and $\delta(h\beta)$ is the discrete delta function.}\\[0.2cm]

Then taking into account (4.7), (4.8) and theorems 4.1, 4.2, for optimal
coefficients we have
$$
C_\beta =D_m(h\beta )*u(h\beta ).\eqno     (4.13)
$$

Thus if we will find the function $u(h\beta )$ then the optimal
coefficients will be found from equality (4.13).

In order to calculate the convolution (4.13) it is required to
find the representation of the function $u(h\beta )$ for all integer
values of $\beta$. From equality (4.1) we get that $u(h\beta ) =
f_m(h\beta )$ when $h\beta \in [0,1]$. Now we need to find
the representation of the function $u(h\beta )$ when $\beta < 0$ and
$\beta>N$.

Since $C_\beta= 0$ when $h\beta \notin [0,1]$ then
$$
C_\beta   = D_m(h\beta )*u(h\beta ) = 0,\,\,\,\,\,\,\,\,\,h\beta
\notin [0,1].
$$

Now we calculate the convolution $v(h\beta ) = G(h\beta)*C_\beta$
when $h\beta \notin [0,1]$.

Suppose $\beta  < 0$ then taking into account equalities (2.9),
(4.2), (4.3), (4.4), we have
$$
v(h\beta ) = G(h\beta )*C_\beta   = \sum\limits_{\gamma  =  -
\infty }^\infty  {C_\gamma  \,G(h\beta  - h\gamma ) = }
$$
$$
= \sum\limits_{\gamma  = 0}^N {C_\gamma  {{{\rm{sign}}(h\beta  -
h\gamma )} \over 2}} \left(\frac{e^{h\beta-h\gamma}-e^{-h\beta+h\gamma}}{2}-
\sum\limits_{k=1}^{m-1}\frac{(h\beta-h\gamma)^{2k-1}}{(2k-1)!}\right) =
$$
$$
=- {1 \over 2}\sum\limits_{\gamma  = 0}^N {C_\gamma  } \left( {{{e^{h\beta  - h\gamma }  - e^{ - h\beta  + h\gamma } } \over 2} - \sum\limits_{k = 1}^{m - 1} {{{\left( {h\beta  - h\gamma } \right)^{2k - 1} } \over {(2k - 1)!}}} } \right) =
$$
$$
=-{{e^{h\beta } } \over 4}\sum\limits_{\gamma  = 0}^N {C_\gamma  e^{ - h\gamma }  + } {{e^{ - h\beta } } \over 4}\sum\limits_{\gamma  = 0}^N {C_\gamma  e^{h\gamma }  + {1 \over 2}\sum\limits_{\gamma  = 0}^N {C_\gamma  } } \sum\limits_{k = 1}^{m - 1} {{{\left( {h\beta  - h\gamma }\right)^{2k - 1} } \over {(2k - 1)!}} = }
$$
$$
= -{{e^{h\beta } } \over 4}(1-e^{-1})  + {{e^{ - h\beta } } \over 4}\sum\limits_{\gamma  = 0}^N {C_\gamma  e^{h\gamma }  + {1 \over 2}\sum\limits_{\gamma  = 0}^N {C_\gamma  } } \left( {\sum\limits_{k = 1}^{\left[ {{{m + 1} \over 2}} \right] - 1} {\sum\limits_{\alpha  = 0}^{2k - 1} {{{(h\beta )^{2k - 1 - \alpha } ( - h\gamma )^\alpha  } \over {(2k - 1 - \alpha )! \cdot \alpha !}} + } } } \right.
$$
$$
 + \left. {\sum\limits_{k = \left[ {{{m + 1} \over 2}} \right]}^{m - 1} {\sum\limits_{\alpha  = 0}^{m - 2} {{{(h\beta )^{2k - 1 - \alpha } ( - h\gamma )^\alpha  } \over {(2k - 1 - \alpha )! \cdot \alpha !}} + } } \sum\limits_{k = \left[ {{{m + 1} \over 2}} \right]}^{m - 1} {\sum\limits_{\alpha  = m - 1}^{2k - 1} {{{(h\beta )^{2k - 1 - \alpha } ( - h\gamma )^\alpha  } \over {(2k - 1 - \alpha )! \cdot \alpha !}}} } } \right) =
$$

$$
 =  - {{e^{h\beta } } \over 4}(1-e^{-1})+{{e^{ - h\beta } } \over 4}\sum\limits_{\gamma  = 0}^N {C_\gamma  e^{h\gamma }  + {1 \over 2}\sum\limits_{k = \left[ {{{m + 1} \over 2}} \right]}^{m - 1} {\sum\limits_{\alpha  = m - 1}^{2k - 1} {{{(h\beta )^{2k - 1 - \alpha } ( - 1)^\alpha  } \over {(2k - 1 - \alpha )! \cdot \alpha !}}} } \sum\limits_{\gamma  = 0}^N {C_\gamma  } } (h\gamma )^\alpha   +
$$
$$
+{1 \over 2}\left({\sum\limits_{k = 1}^{
\left[ {{{m + 1} \over 2}} \right]-1}{\sum\limits_{\alpha  = 0}^{2k - 1}
{{{(h\beta )^{2k - 1 - \alpha } ( - 1)^\alpha  }\over {(2k - 1 - \alpha )!
\cdot (\alpha+1)!}}+\sum\limits_{k=\left[{{{m + 1} \over 2}}\right]}^{m - 1}
{\sum\limits_{\alpha=0}^{m - 2} {{{(h\beta )^{2k - 1 - \alpha }(-1)^\alpha }
\over {(2k - 1 - \alpha )! \cdot (\alpha+1)!}}} } } } } \right) =
$$
$$
=-{{e^{h\beta } } \over 4}(1-e^{-1})+De^{-h\beta }+Q^{(2m - 3)}(h\beta ) + Q_{m - 2} (h\beta ).
$$
Thus we get
$$
v(h\beta)=-{{e^{h\beta } } \over 4}
(1-e^{-1})+De^{-h\beta }+Q^{(2m - 3)}(h\beta ) + Q_{m - 2} (h\beta ).
\eqno (4.14)
$$
where
$$
Q^{(2m - 3)}(h\beta )={1 \over 2}
\left({\sum\limits_{k = 1}^{\left[ {{{m + 1} \over 2}} \right] - 1}
{\sum\limits_{\alpha  = 0}^{2k - 1} {{{(h\beta )^{2k - 1 - \alpha }
(-1)^\alpha  } \over {(2k - 1 - \alpha )! \cdot (\alpha+1)!}}   + } } } \right.
$$
$$
+\left. {\sum\limits_{k = \left[ {{{m + 1} \over 2}} \right]}^{m - 1}
 {\sum\limits_{\alpha  = 0}^{m - 2} {{{(h\beta )^{2k - 1 - \alpha }
 ( - 1)^\alpha} \over {(2k - 1 - \alpha )! \cdot (\alpha+1)!}}} } } \right)
\eqno (4.15)
$$
is the polynomial of degree $2m-3$ with respect to $(h\beta)$,
$$
Q_{m - 2}(h\beta )={1 \over 2}\sum\limits_{k=
\left[{{{m + 1} \over 2}} \right]}^{m - 1}
{\sum\limits_{\alpha  = m - 1}^{2k - 1} {{{(h\beta )^{2k - 1 - \alpha }
(-1)^\alpha  } \over {(2k - 1 - \alpha )! \cdot \alpha !}}} }
\sum\limits_{\gamma  = 0}^N {C_\gamma  } (h\gamma )^\alpha\eqno (4.16)
$$
is unknown polynomial of degree $m-2$ with respect to $(h\beta)$,
$$
D={1 \over 4}\sum\limits_{\gamma  = 0}^N {C_\gamma  e^{h\gamma } }.\eqno (4.17)
$$

Similarly, in the case $\beta>N$ for the convolution $v(h\beta ) = G(h\beta)*C_\beta$
we obtain
$$
v(h\beta )= {{e^{h\beta } } \over 4}(1-e^{-1})-De^{ - h\beta }-
Q^{(2m - 3)} (h\beta ) - Q_{m - 2} (h\beta ).\eqno (4.18)
$$

We denote
$$
Q_{m - 2}^{(-)}(h\beta )=
P_{m - 2} (h\beta ) + Q_{m - 2} (h\beta ),\,\,\,\,\,\,a^ -   = d + D,\eqno (4.19)
$$

    $$
Q_{m - 2}^{( + )} (h\beta ) = P_{m - 2} (h\beta ) -
Q_{m - 2} (h\beta ),\,\,\,\,\,\,a^ +   = d - D.\eqno (4.20)
$$
and taking into account (4.14), (4.18), (4.7) we get the following
problem.

\textbf{Problem B.} \emph{Find the solution of the equation}
$$
D_m(h\beta )*u(h\beta ) = 0,\,\,\,\,\,\,\,\,h\beta  \notin [0,1]\eqno
(4.21)
$$
\emph{having the form:}
$$
u(h\beta ) =
\left\{
\begin{array}{ll}
-\frac{e^{h\beta}}{4}(1-e^{-1})+a^-e^{-h\beta}+
Q^{(2m-3)}(h\beta)+Q_{m-2}^{(-)}(h\beta),& \beta<0; \\
f_m(h\beta ),& 0 \le \beta  \le N; \\
\frac{e^{h\beta}}{4}(1-e^{-1})+a^+e^{-h\beta}-
Q^{(2m-3)}(h\beta)+Q_{m-2}^{(+)}(h\beta),& \beta>N.\\
\end{array}
\right.\eqno (4.22)
$$
\emph{Here $Q_{m-2}^{(-)}(h\beta)$ and $Q_{m-2}^{(+)}(h\beta)$ are
unknown polynomials of degree $m-2$ with respect to $h\beta$,
$a^-$ and $a^+$ are unknown constants.}\\[0.2cm]

If we find $Q_{m-2}^{(-)}(h\beta)$, $Q_{m-2}^{(+)}(h\beta)$, $a^-$
and $a^+$ then from (4.19), (4.20) we have
$$
P_{m-2}(h\beta)=\frac{1}{2}\left(Q_{m-2}^{(-)}(h\beta)+Q_{m-2}^{(+)}(h\beta)\right),\ \
d=\frac{1}{2}(a^-+a^+),
$$
$$
Q_{m-2}(h\beta)=\frac{1}{2}\left(Q_{m-2}^{(-)}(h\beta)-Q_{m-2}^{(+)}(h\beta)\right),\
\ D=\frac{1}{2}(a^--a^+),
$$

Unknowns $Q_{m-2}^{(-)}(h\beta)$, $Q_{m-2}^{(+)}(h\beta)$, $a^-$ and $a^+$
can be found  from the equation (4.21), using the function
$D_m(h\beta)$. Then we can obtain explicit form of the function
$u(h\beta )$ and find the optimal coefficients $C_\beta$.
Thus, the problem $B$ and respectively the problem $A$ can be solved.

But here we will not find $Q_{m-2}^{(-)}(h\beta)$,
$Q_{m-2}^{(+)}(h\beta)$, $a^-$ and $a^+$. Instead of them, using $D_m(h\beta)$ and
$u(h\beta)$, taking into account (4.13), we will find the expressions
for optimal coefficients $C_{\beta}$ when $\beta=1,...,N-1$.

We denote
$$
a_k={{A_k}\over {\lambda _k p}}\sum\limits_{\gamma  = 1}^\infty
{\lambda _k^\gamma  } \left( { - {{e^{ - h\gamma } } \over 4}
(1-e^{-1})  + Q^{(2m - 3)} ( - h\gamma ) + Q_{m - 2}^{( - )} ( - h\gamma ) + a^ -  e^{h\gamma }  -
f_m ( - h\gamma )} \right),\eqno (4.23)
$$
$$
b_k={{A_k }\over{\lambda _k p}}\sum\limits_{\gamma  = 1}^\infty
{\lambda _k^\gamma  } \left( {{{e^{h\gamma  + 1} } \over 4}
(1-e^{-1})- Q^{(2m - 3)} (1 + h\gamma )
+ Q_{m - 2}^{( + )} (1 + h\gamma ) + a^ +  e^{ - 1 - h\gamma }  -
f_m (1 + h\gamma )} \right),\eqno (4.24)
$$
Here $\lambda_k$ are the roots  and $p$ is the leading coefficient
of the polynomial $\mathcal{P}_{2m-2}(\lambda)$ of degree $2m-2$
defined by equality (4.12) and $|\lambda_k|<1$. The series in the
notations (4.23), (4.24) are convergent.

The following is true

\textbf{Theorem 4.3} (Theorem 3 of [23]). {\it The coefficients of
optimal quadrature formulas in the sense of Sard of the form (1.1)
in the space $W_2^{(m,m-1)}(0,1)$ have the following form
$$
C_\beta=D_m(h\beta)*f_m(h\beta)+\sum\limits_{k=1}^{m-1}
\left(a_k\lambda_k^{\beta}+b_k\lambda_k^{N-\beta}\right)
,\ \ \beta=1,2,...,N-1, \eqno (4.25)
$$
where $a_k$ and $b_k$ are unknowns and have the form (4.23) and
(4.24) respectively, $\lambda_k$ are the roots of the polynomial
$\mathcal{P}_{2m-2}(\lambda)$ which defined by equality (4.12) and
$|\lambda_k|<1$.}

{\bf Proof.} Suppose $\beta=\overline{1,N-1}$. Then from (4.13),
using (4.9), (4.22) and theorem 4.2, we have
$$
C_{\beta}=D_m(h\beta)*u(h\beta)=\sum\limits_{\gamma=-\infty}^{\infty}
D_m(h\beta-h\gamma)u(h\gamma)=
$$
$$
=\sum\limits_{\gamma=-\infty}^{-1}D_m(h\beta-h\gamma)u(h\gamma)+
\sum\limits_{\gamma=0}^ND_m(h\beta-h\gamma)u(h\gamma)+
\sum\limits_{\gamma=N+1}^{\infty}D_m(h\beta-h\gamma)u(h\gamma).
$$
Hence using the definition of convolution of discrete functions we get
$$
C_{\beta}=D_m(h\beta)*f_m(h\beta)+
$$
$$
+\sum\limits_{\gamma=1}^\infty{\sum\limits_{k = 1}^{m - 1} {{{A_k }
\over {\lambda _k p}}\lambda _k^{\beta  + \gamma } } \left( { - {{e^{ - h\gamma } }
\over 4}(1-e^{-1})+Q^{(2m - 3)} ( - h\gamma ) +
Q_{m - 2}^{( - )} ( - h\gamma ) + a^ -  e^{h\gamma }  - f_m ( - h\gamma )} \right)}  +
$$
$$
+\sum\limits_{\gamma  = 1}^\infty
\sum\limits_{k = 1}^{m - 1} {{{A_k } \over {\lambda _k p}}
\lambda _k^{N + \gamma  - \beta } } \Bigg( {{e^{1 + h\gamma } }
\over 4}(1-e^{-1})-Q^{(2m - 3)}(1+h\gamma) + Q_{m - 2}^{( + )}
(1 + h\gamma ) + a^ +  e^{ - 1 - h\gamma }  -
$$
$$
-f_m (1 + h\gamma ) \Bigg).
$$
Whence taking into account the notations (4.23), (4.24), we get (4.25).

Theorem 4.3 is proved.\hfill $\Box$ \\[0.05 cm]

From theorem 4.3 it is clear that in order to obtain explicit forms of the
optimal coefficients $C_{\beta}$ in the space $W_2^{(m.m-1)}(0,1)$ it is
sufficient to find $a_k$ and $b_k$ ($k=\overline{1,m-1}$). But here we will not to calculate
series (4.23) and (4.24). Instead of them substituting the
equality (4.25) into (4.1) we obtain identity with respect to
$(h\beta)$. Whence, equating corresponding coefficients the left
and the right hand sides of equation (4.1) we will find $a_k$
and $b_k$. And the coefficient $C_0$ and $C_N$ will be found from
(4.3) when $\alpha=0$ and (4.4), respectively. Below we will do it.

It should be noted that the cases $m=1$ and $m=2$ are solved in
the work [23] and the following theorems are proved.

\textbf{Theorem 4.4} (Theorem 4 of [23]). {\it The coefficients of
optimal quadrature formulas of the form (1.1) with equal spaced
nodes in the space $W_2^{(1,0)}(0,1)$ are expressed by formulas
$$
C_{\beta}=\left\{
\begin{array}{ll}
\frac{e^h-1}{e^h+1},& \beta=0,N,\\
\frac{2(e^h-1)}{e^h+1},& \beta=\overline{1,N-1},\\
\end{array}
\right.
$$
where $h=1/N$, $N=1,2,...$. }

\textbf{Theorem 4.5} (Theorem 5 of [23]). {\it The coefficients of
optimal quadrature formulas of the form (1.1) with equal spaced
nodes in the space $W_2^{(2,1)}(0,1)$ are expressed by formulas
$$
C_{\beta}=\left\{
\begin{array}{ll}
1-\frac{h}{e^h-1}-K(h)(\lambda_1-\lambda_1^N),& \beta=0,\\
h+K(h)\left((e^h-\lambda_1)\lambda_1^{\beta}+(1-\lambda_1e^h)\lambda_1^{N-\beta}\right),
& \beta=\overline{1,N-1},\\
-1+\frac{e^hh}{e^h-1}-K(h)(\lambda_1-\lambda_1^N)e^h,& \beta=N,\\
\end{array}
\right.
$$
where
$$
K(h)=\frac{(2e^h-2-he^h-h)(\lambda_1-1)}{2(e^h-1)^2(\lambda_1+\lambda_1^{N+1})},
$$
$$
\lambda_1=\frac{h(e^{2h}+1)-e^{2h}+1-(e^h-1)\sqrt{h^2(e^{h}+1)^2+2h(1-e^h)}}
{1-e^{2h}+2he^h},\ \ |\lambda_1|<1,
$$
$h=1/N$, $N=2,3,...$.}\\[0.05 cm]

The main goal of the present section is to solve the system
(4.1)-(4.4) for any $m\geq 2$ and any natural $N$, $N\geq m$. As
mentioned above for this sufficiently to find $a_k$ and $b_k$
($k=\overline{1,m-1}$) in (4.25).

\textbf{The main result} of the present paper is the following theorem.

\textbf{Theorem 4.6.} {\it The coefficients of optimal quadrature formulas of the form (1.1)
with the error functional (1.2) and with equal spaced nodes in the space $W_2^{(m,m-1)}(0,1)$
when $m\geq 2$ are expressed by  formulas}\\
$ C_0= \frac{e^h-1-h}{e^h-1}+\sum\limits_{k=1}^{m-1}
\Bigg(a_k\frac{\lambda_k(e^h-e)+\lambda_k^2(e-1)+\lambda_k^{N+1}(1-e^h)}
{(e-1)(1-\lambda_k)(e^h-\lambda_k)}+b_k\frac{\lambda_k^{N+1}(e^h-e)+\lambda_k^N(e-1)+\lambda_k(1-e^h)}
{(e-1)(\lambda_k-1)(\lambda_ke^h-1)}\Bigg),
$ \\
$
C_\beta=
h+\sum\limits_{k=1}^{m-1}\left(a_k\lambda_k^\beta+b_k\lambda_k^{N-\beta}\right),
\ \ \ \ \beta=\overline{1,N-1},\\
$\\
$
C_N=
\frac{e^hh+1-e^h}{e^h-1}+\sum\limits_{k=1}^{m-1}
\Bigg(a_k\frac{\lambda_k(e-e^{h+1})+\lambda_k^N(e^{h+1}-e^h)+\lambda_k^{N+1}(e^h-e)}
{(e-1)(1-\lambda_k)(e^h-\lambda_k)}+
b_k\frac{\lambda_k^{N+1}(e-e^{h+1})+\lambda_k^2(e^{h+1}-e^h)+\lambda_k(e^h-e)}
{(e-1)(1-\lambda_k)(1-\lambda_ke^h)}\Bigg),
$ \\
\emph{where $a_k$ and $b_k$ ($k=\overline{1,m-1}$) are defined by the following system of $2m-2$ linear equations}\\
$$
\begin{array}{l}
\sum\limits_{k=1}^{m-1}a_k\frac{\lambda_k}{(\lambda_k-1)(\lambda_k-e^h)}+\sum\limits_{k=1}^{m-1}b_k\frac{\lambda_k^{N+1}}{(\lambda_k-1)
(\lambda_ke^h-1)}=\frac{h-2}{2(e^h-1)}+\frac{h}{(e^h-1)^2};\\
\sum\limits_{k=1}^{m-1}a_k\frac{\lambda_k^{N+1}}{(\lambda_k-1)(\lambda_k-e^h)}+\sum\limits_{k=1}^{m-1}b_k\frac{\lambda_k}{(\lambda_k-1)
(\lambda_ke^h-1)}=\frac{h-2}{2(e^h-1)}+\frac{h}{(e^h-1)^2};\\
\sum\limits_{k=1}^{m-1}a_k\Bigg[\sum\limits_{l=2}^{j}\frac{h^{2l-2}}{(2l-2)!}
\sum\limits_{i=1}^{2l-2}\frac{\lambda_k\Delta^i0^{2l-2}}{(\lambda_k-1)^{i+1}}\Bigg]+
\sum\limits_{k=1}^{m-1}b_k\Bigg[\sum\limits_{l=2}^{j}\frac{h^{2l-2}}{(2l-2)!}
\sum\limits_{i=1}^{2l-2}\frac{\lambda_k^{N+i}\Delta^i0^{2l-2}}{(1-\lambda_k)^{i+1}}\Bigg]=0,\
\ j=\overline{2,[\frac{m}{2}]};\\
\sum\limits_{k=1}^{m-1}a_k\Bigg[\sum\limits_{l=1}^{j}\frac{h^{2l-1}}{(2l-1)!}
\sum\limits_{i=1}^{2l-1}\frac{\lambda_k\Delta^i0^{2l-1}}{(\lambda_k-1)^{i+1}}\Bigg]+
\sum\limits_{k=1}^{m-1}b_k\Bigg[\sum\limits_{l=1}^{j}\frac{h^{2l-1}}{(2l-1)!}
\sum\limits_{i=1}^{2l-1}\frac{\lambda_k^{N+i}\Delta^i0^{2l-1}}{(1-\lambda_k)^{i+1}}\Bigg]=\sum\limits_{l=1}^j\frac{h^{2l}B_{2l}}{(2l)!},\
\ j=\overline{1,[\frac{m-1}{2}]};\\
\end{array}
$$
$$
\begin{array}{l}
\sum\limits_{k=1}^{m-1}a_k\Bigg[\sum\limits_{l=1}^{j}h^lC_j^l
\sum\limits_{i=1}^{l}\frac{\lambda_k^{N+i}\Delta^i0^{l}}{(1-\lambda_k)^{i+1}}-h^j
\sum\limits_{i=1}^j\frac{\lambda_k^i\Delta^i0^j}{(1-\lambda_k)^{i+1}}\Bigg]\\
\ \ \ \ \ \ \ \
+\sum\limits_{k=1}^{m-1}b_k\Bigg[\sum\limits_{l=1}^{j}h^lC_j^l
\sum\limits_{i=1}^{l}\frac{\lambda_k\Delta^i0^{l}}{(\lambda_k-1)^{i+1}}-h^j
\sum\limits_{i=1}^j\frac{\lambda_k^{N+1}\Delta^i0^j}{(\lambda_k-1)^{i+1}}\Bigg]=\sum\limits_{l=1}^{j-1}\frac{j!B_{j+1-l}}{l!(j+1-l)!}h^{j+1-l},\
\ j=\overline{1,m-2}.
\\
\end{array}
$$
\emph{Here $\lambda_k$ are the roots of the polynomial (4.12) and
$|\lambda_k|<1$, $B_j$ are Bernoulli
numbers.} \\

In the proof of theorem 4.6 we use the following formulas from
[10]:
$$
\sum\limits_{\gamma  = 0}^{n - 1} {q^\gamma  \gamma ^k  =
\frac{1}{{1 - q}}\sum\limits_{i = 0}^k {\left( {\frac{q}{{1 - q}}}
\right)^i \Delta ^i 0^k  - \frac{{q^n }}{{1 - q}}\sum\limits_{i =
0}^k {\left( {\frac{q}{{1 - q}}} \right)^i \Delta ^i \gamma ^k
|_{\gamma  = n} ,} } }.\eqno    (4.26)
$$
where  $\Delta ^i 0^k   = \sum\limits_{l =
1}^i {( - 1)^{i - l} C_i^l l^k }$, $\Delta^i\gamma^k$ is the finite difference of order $i$ of
$\gamma^k$,\\
and from [9]:
$$
\sum\limits_{\gamma  = 0}^{\beta  - 1} {\gamma ^k  =
\sum\limits_{j = 1}^{k + 1} {\frac{{k!\,B_{k + 1 - j} }}{{j!\,(k +
1 - j)!}}\,\beta ^j ,} }\eqno    (4.27)
$$
where $B_{k + 1 - j} $ are Bernoulli numbers,
    $$
\Delta ^\alpha  x^\nu   = \sum\limits_{p = 0}^\nu  {\left(
{\begin{array}{c}
   \nu   \\
   p  \\
\end{array}} \right)\Delta ^\alpha  } 0^p x^{\nu  - p}.
\eqno (4.28)
$$

\textbf{Proof of theorem 4.6}. For convenience we denote
$$
T=D_m(h\beta)*f_m(h\beta).\eqno (4.29)
$$
Now we consider equality (4.1)
$$
\sum\limits_{\gamma=0}^NC_{\gamma}\frac{\mathrm{sign}(h\beta-h\gamma)}{2}
\left(\frac{e^{h\beta-h\gamma}-e^{h\gamma-h\beta}}{2}-\sum\limits_{k=1}^{m-1}
\frac{(h\beta-h\gamma)^{2k-1}}{(2k-1)!}\right)+P_{m-2}(h\beta)+de^{-h\beta}=f_m(h\beta),
\eqno (4.30)
$$
where $\beta=0,1,...,N$.

We denote
$$
g(h\beta)=\sum\limits_{\gamma=0}^NC_{\gamma}\frac{\mathrm{sign}(h\beta-h\gamma)}{2}
\left(\frac{e^{h\beta-h\gamma}-e^{h\gamma-h\beta}}{2}-\sum\limits_{k=1}^{m-1}
\frac{(h\beta-h\gamma)^{2k-1}}{(2k-1)!}\right).\eqno (4.31)
$$
The expresson $g(h\beta)$ we rewrite in the following form
$$
\begin{array}{ll}
g(h\beta)=&\displaystyle C_0\left(\frac{e^{h\beta}-e^{-h\beta}}{2}-
\sum\limits_{k=1}^{m-1}\frac{(h\beta)^{2k-1}}{(2k-1)!}\right)+
\sum\limits_{\gamma=1}^{\beta-1}C_{\gamma}\left(\frac{e^{h\beta-h\gamma}-e^{h\gamma-h\beta}}{2}-\sum\limits_{k=1}^{m-1}
\frac{(h\beta-h\gamma)^{2k-1}}{(2k-1)!}\right) \\
 &\displaystyle -\frac{1}{2}
\sum\limits_{\gamma=0}^{N}C_{\gamma}\left(\frac{e^{h\beta-h\gamma}-e^{h\gamma-h\beta}}{2}-\sum\limits_{k=1}^{m-1}
\frac{(h\beta-h\gamma)^{2k-1}}{(2k-1)!}\right).
\end{array}\eqno (4.32)
$$
Further we denote
$$
g_1(h\beta)=\sum\limits_{\gamma=0}^{\beta-1}C_{\gamma}\left(\frac{e^{h\beta-h\gamma}-e^{h\gamma-h\beta}}{2}-\sum\limits_{k=1}^{m-1}
\frac{(h\beta-h\gamma)^{2k-1}}{(2k-1)!}\right),\eqno (4.33)
$$
$$
g_2(h\beta)=-\frac{1}{2}\sum\limits_{\gamma=0}^{N}C_{\gamma}\left(\frac{e^{h\beta-h\gamma}-e^{h\gamma-h\beta}}{2}-\sum\limits_{k=1}^{m-1}
\frac{(h\beta-h\gamma)^{2k-1}}{(2k-1)!}\right).\eqno (4.34)
$$
Firstly we consider $g_1(h\beta)$ and rewrite (4.33) in the
following form
$$
g_1(h\beta)=\sum\limits_{\gamma=1}^{\beta-1}C_{\beta-\gamma}
\left(\frac{e^{h\gamma}-e^{-h\gamma}}{2}-
\sum\limits_{k=1}^{m-1}\frac{(h\gamma)^{2k-1}}{(2k-1)!}\right).
$$
Hence using (4.25), (4.29) we get
$$
g_1(h\beta)=\sum\limits_{\gamma=1}^{\beta-1}\left(T+
\sum\limits_{k=1}^{m-1}\left(a_k\lambda_k^{\beta-\gamma}+
b_k\lambda_k^{N-\beta+\gamma}\right)\right)\frac{e^{h\gamma}}{2}-
$$
$$
-\sum\limits_{\gamma=1}^{\beta-1}\left(T+
\sum\limits_{k=1}^{m-1}\left(a_k\lambda_k^{\beta-\gamma}+
b_k\lambda_k^{N-\beta+\gamma}\right)\right)\frac{e^{-h\gamma}}{2}-
$$
$$
-\sum\limits_{\ell=1}^{m-1}\sum\limits_{\gamma=1}^{\beta-1}\left(T+
\sum\limits_{k=1}^{m-1}\left(a_k\lambda_k^{\beta-\gamma}+
b_k\lambda_k^{N-\beta+\gamma}\right)\right)\frac{(h\gamma)^{2\ell-1}}{(2\ell-1)!}.
$$
From here taking into account (4.26), (4.27) and after some
simplifications we have
$$
g_1(h\beta)=\frac{T(e^h-e^{h\beta})}{2(1-e^h)}+
\sum\limits_{k=1}^{m-1}\left(a_k\frac{\lambda_k^{\beta}e^h-
\lambda_ke^{h\beta}}{2(\lambda_k-e^h)}+b_k\frac{\lambda_k^{N-\beta+1}e^h-
\lambda_k^Ne^{h\beta}}{2(1-\lambda_ke^h)}\right)-
$$
$$
-\frac{T(1-e^{h-h\beta})}{2(e^h-1)}-
\sum\limits_{k=1}^{m-1}\left(a_k\frac{\lambda_k^{\beta}-
\lambda_ke^{h-h\beta}}{2(\lambda_ke^h-1)}+b_k\frac{\lambda_k^{N-\beta+1}-
\lambda_k^Ne^{h-h\beta}}{2(e^h-\lambda_k)}\right)-
$$
$$
-\sum\limits_{\ell=1}^{m-1}\frac{Th^{-1}(h\beta)^{2l}}{(2l)!}-\sum\limits_{\ell=1}^{m-1}
\sum\limits_{j=1}^{2\ell-1}\frac{Th^{2\ell-1}B_{2\ell-j}}{j!(2\ell-j)!}\beta^j-
$$
$$
-\sum\limits_{\ell=1}^{m-1}\frac{h^{2\ell-1}}{(2\ell-1)!}\sum\limits_{k=1}^{m-1}
\Bigg[a_k\left(\frac{\lambda_k^{\beta+1}}{\lambda_k-1}\sum\limits_{i=0}^{2\ell-1}
\frac{\Delta^i0^{2\ell-1}}{(\lambda_k-1)^i}-\frac{\lambda_k}{\lambda_k-1}
\sum\limits_{i=0}^{2\ell-1}\frac{\Delta^i\beta^{2\ell-1}}{(\lambda_k-1)^i}\right)+
$$
$$
+b_k\left(\frac{\lambda_k^{N-\beta}}{1-\lambda_k}\sum\limits_{i=0}^{2\ell-1}
\left(\frac{\lambda_k}{1-\lambda_k}\right)^i\Delta^i0^{2\ell-1}-
\frac{\lambda_k^N}{1-\lambda_k}
\sum\limits_{i=0}^{2\ell-1}\left(\frac{\lambda_k}
{1-\lambda_k}\right)^i\Delta^i\beta^{2\ell-1}
\right)\Bigg].\eqno (4.35)
$$
In the last expression of $g_1(h\beta)$ the coefficients of
$\lambda_k^{\beta}$ and $\lambda_k^{N-\beta}$ are the values of
the polynomial $\mathcal{P}_{2m-2}(\lambda)$ which defined by
equality (4.12) at the $\lambda_k$. Since $\lambda_k$ are the
roots of the polynomial (4.12), then the coefficients of
$\lambda_k^{\beta}$ and $\lambda_k^{N-\beta}$ are zero. Then from
(4.35) for $g_1(h\beta)$ we get
$$
g_1(h\beta)=\frac{e^{h\beta}}{2}\left[\frac{T}{e^h-1}+\sum\limits_{k=1}^{m-1}
\left(a_k\frac{\lambda_k}{e^h-\lambda_k}+b_k\frac{\lambda_k^N}{\lambda_ke^h-1}\right)\right]-
$$
$$
-\frac{e^{-h\beta}}{2}\left[\frac{Te^h}{1-e^h}+\sum\limits_{k=1}^{m-1}
\left(a_k\frac{\lambda_ke^h}{1-e^h\lambda_k}+b_k\frac{\lambda_k^Ne^h}{\lambda_k-e^h}\right)
\right]+
$$
$$
+\frac{Te^h}{2(1-e^h)}+\frac{T}{2(1-e^h)}-\sum\limits_{\ell=1}^{m-1}
\frac{Th^{-1}(h\beta)^{2\ell}}{(2\ell)!}-\sum\limits_{\ell=1}^{m-1}Th^{2\ell-1}
\sum\limits_{j=1}^{2\ell-1}\frac{B_{2\ell-j}}{j!(2\ell-j)!}\beta^j+
$$
$$
+\sum\limits_{\ell=1}^{m-1}\frac{h^{2\ell-1}}{(2\ell-1)!}\sum\limits_{k=1}^{m-1}
a_k\frac{\lambda_k}{\lambda_k-1}\sum\limits_{i=0}^{2\ell-1}\frac{\Delta^i\beta^{2\ell-1}}
{(\lambda_k-1)^i}+
\sum\limits_{\ell=1}^{m-1}\frac{h^{2\ell-1}}{(2\ell-1)!}\sum\limits_{k=1}^{m-1}
b_k\frac{\lambda_k^N}{1-\lambda_k}\sum\limits_{i=0}^{2\ell-1}
\left(\frac{\lambda_k}{1-\lambda_k}\right)^i\Delta^i\beta^{2\ell-1}.
$$
Finally taking into account equality (4.28) we have
$$
g_1(h\beta)=\frac{e^{h\beta}}{2}\left[\frac{T}{e^h-1}+\sum\limits_{k=1}^{m-1}
\left(a_k\frac{\lambda_k}{e^h-\lambda_k}+b_k\frac{\lambda_k^N}{\lambda_ke^h-1}\right)\right]-
$$
$$
-\frac{e^{-h\beta}}{2}\left[\frac{Te^h}{1-e^h}+\sum\limits_{k=1}^{m-1}
\left(a_k\frac{\lambda_ke^h}{1-e^h\lambda_k}+b_k\frac{\lambda_k^Ne^h}{\lambda_k-e^h}\right)
\right]+
$$
$$
+\frac{Te^h}{2(1-e^h)}+\frac{T}{2(1-e^h)}-\sum\limits_{\ell=1}^{m-1}
\frac{Th^{-1}(h\beta)^{2\ell}}{(2\ell)!}-\sum\limits_{\ell=1}^{m-1}Th^{2\ell-1}
\sum\limits_{j=1}^{2\ell-1}\frac{B_{2\ell-j}}{j!(2\ell-j)!}\beta^j+
$$
$$
+\sum\limits_{\ell=1}^{m-1}\frac{h^{2\ell-1}}{(2\ell-1)!}
\sum\limits_{j=0}^{2\ell-1}C_{2\ell-1}^j\beta^j
\sum\limits_{k=1}^{m-1}
a_k\frac{\lambda_k}{\lambda_k-1}\sum\limits_{i=0}^{2\ell-1}\frac{\Delta^i0^{2\ell-1-j}}
{(\lambda_k-1)^i}+
$$
$$
+\sum\limits_{\ell=1}^{m-1}\frac{h^{2\ell-1}}{(2\ell-1)!}
\sum\limits_{j=0}^{2\ell-1}C_{2\ell-1}^j\beta^j
\sum\limits_{k=1}^{m-1}
b_k\frac{\lambda_k^N}{1-\lambda_k}\sum\limits_{i=0}^{2\ell-1}
\left(\frac{\lambda_k}{1-\lambda_k}\right)^i\Delta^i0^{2\ell-1-j}.\eqno (4.36)
$$
Now using the binomial formula and equalities (4.3), (4.4) from (4.34) we obtain
$$
g_2(h\beta)=-\frac{1}{2}\left(\frac{e^{h\beta}}{2}
\sum\limits_{\gamma=0}^NC_{\gamma}e^{-h\gamma}-\frac{e^{-h\beta}}{2}
\sum\limits_{\gamma=0}^NC_{\gamma}e^{h\gamma}-
\sum\limits_{\gamma=0}^NC_{\gamma}\sum\limits_{k=1}^{m-1}
\frac{(h\beta-h\gamma)^{2m-1}}{(2k-1)!}\right)=
$$
$$
=-\frac{1}{2}\Bigg[\frac{e^{h\beta}}{2}(1-e^{-1})-
\frac{e^{-h\beta}}{2}\sum\limits_{\gamma=0}^NC_{\gamma}e^{h\gamma}
-\Bigg(\sum\limits_{k=1}^{\left[\frac{m+1}{2}\right]-1}
\sum\limits_{\alpha=0}^{2k-1}\frac{(h\beta)^{2k-1-\alpha}(-1)^{\alpha}}
{(2k-1-\alpha)!(\alpha+1)!}+
$$
$$
+\sum\limits_{k=\left[\frac{m+1}{2}\right]}^{m-1}
\sum\limits_{\alpha=0}^{m-2}\frac{(h\beta)^{2k-1-\alpha}(-1)^{\alpha}}
{(2k-1-\alpha)!(\alpha+1)!}+
\sum\limits_{k=\left[\frac{m+1}{2}\right]}^{m-1}
\sum\limits_{\alpha=m-1}^{2k-1}\frac{(h\beta)^{2k-1-\alpha}(-1)^{\alpha}}
{(2k-1-\alpha)!\alpha!}\sum\limits_{\gamma=0}^NC_{\gamma}(h\gamma)^{\alpha}\Bigg)\Bigg].
\eqno (4.37)
$$

Using the binomial formula in equality (4.5) for $f_m(h\beta)$ we have
$$
f_m(h\beta)=\frac{e^{h\beta}+e^{-h\beta}+e^{1-h\beta}+e^{h\beta-1}-4}
{4}-\sum\limits_{k=1}^{m-1}\frac{(h\beta)^{2k}}{(2k)!}+
\frac{1}{2}\sum\limits_{k=1}^{\left[\frac{m+1}{2}\right]-1}\sum\limits_{\alpha=0}^{2k-1}
\frac{(h\beta)^{2k-1-\alpha}(-1)^{\alpha}}{(2k-1-\alpha)!(\alpha+1)!}+
$$
$$
+\frac{1}{2}\sum\limits_{k=\left[\frac{m+1}{2}\right]}^{m-1}
\sum\limits_{\alpha=0}^{m-2}
\frac{(h\beta)^{2k-1-\alpha}(-1)^{\alpha}}{(2k-1-\alpha)!(\alpha+1)!}+
\frac{1}{2}\sum\limits_{k=\left[\frac{m+1}{2}\right]}^{m-1}
\sum\limits_{\alpha=m-1}^{m-2}
\frac{(h\beta)^{2k-1-\alpha}(-1)^{\alpha}}{(2k-1-\alpha)!(\alpha+1)!}.\eqno (4.38)
$$

Taking into account (4.36), (4.37) and putting (4.32), (4.38) into (4.30) we get
$$
C_0\left[\frac{e^{h\beta}-e^{-h\beta}}{2}-\sum\limits_{k=1}^{m-1}
\frac{(h\beta)^{2k-1}}{(2k-1)!}\right]+
\frac{e^{h\beta}}{2}\left[\frac{T}{e^h-1}+\sum\limits_{k=1}^{m-1}
\left(a_k\frac{\lambda_k}{e^h-\lambda_k}+b_k\frac{\lambda_k^N}{\lambda_ke^h-1}\right)\right]-
$$
$$
-\frac{e^{-h\beta}}{2}\left[\frac{Te^h}{1-e^h}+\sum\limits_{k=1}^{m-1}
\left(a_k\frac{\lambda_ke^h}{1-e^h\lambda_k}+b_k\frac{\lambda_k^Ne^h}{\lambda_k-e^h}\right)
\right]+
$$
$$
+\frac{Te^h}{2(1-e^h)}+\frac{T}{2(1-e^h)}-\sum\limits_{\ell=1}^{m-1}
\frac{Th^{-1}(h\beta)^{2\ell}}{(2\ell)!}-\sum\limits_{\ell=1}^{m-1}Th^{2\ell-1}
\sum\limits_{j=1}^{2\ell-1}\frac{B_{2\ell-j}}{j!(2\ell-j)!}\beta^j+
$$
$$
+\sum\limits_{\ell=1}^{m-1}\frac{h^{2\ell-1}}{(2\ell-1)!}
\sum\limits_{j=0}^{2\ell-1}C_{2\ell-1}^j\beta^j
\sum\limits_{k=1}^{m-1}
a_k\frac{\lambda_k}{\lambda_k-1}\sum\limits_{i=0}^{2\ell-1}\frac{\Delta^i0^{2\ell-1-j}}
{(\lambda_k-1)^i}+
$$
$$
+\sum\limits_{\ell=1}^{m-1}\frac{h^{2\ell-1}}{(2\ell-1)!}
\sum\limits_{j=0}^{2\ell-1}C_{2\ell-1}^j\beta^j
\sum\limits_{k=1}^{m-1}
b_k\frac{\lambda_k^N}{1-\lambda_k}\sum\limits_{i=0}^{2\ell-1}
\left(\frac{\lambda_k}{1-\lambda_k}\right)^i\Delta^i0^{2\ell-1-j}-
$$
$$
-\frac{1}{2}\Bigg[\frac{e^{h\beta}}{2}(1-e^{-1})-
\frac{e^{-h\beta}}{2}\sum\limits_{\gamma=0}^NC_{\gamma}e^{h\gamma}
-\Bigg(\sum\limits_{k=1}^{\left[\frac{m+1}{2}\right]-1}
\sum\limits_{\alpha=0}^{2k-1}\frac{(h\beta)^{2k-1-\alpha}(-1)^{\alpha}}
{(2k-1-\alpha)!(\alpha+1)!}+
$$
$$
+\sum\limits_{k=\left[\frac{m+1}{2}\right]}^{m-1}
\sum\limits_{\alpha=0}^{m-2}\frac{(h\beta)^{2k-1-\alpha}(-1)^{\alpha}}
{(2k-1-\alpha)!(\alpha+1)!}+
\sum\limits_{k=\left[\frac{m+1}{2}\right]}^{m-1}
\sum\limits_{\alpha=m-1}^{2k-1}\frac{(h\beta)^{2k-1-\alpha}(-1)^{\alpha}}
{(2k-1-\alpha)!\alpha!}\sum\limits_{\gamma=0}^NC_{\gamma}(h\gamma)^{\alpha}\Bigg)\Bigg]+
$$
$$
+P_{m-2}(h\beta)+2d\frac{e^{-h\beta}}{2}=\frac{e^{h\beta}+e^{-h\beta}+e^{1-h\beta}+e^{h\beta-1}-4}
{4}-\sum\limits_{k=1}^{m-1}\frac{(h\beta)^{2k}}{(2k)!}+
$$
$$
+\frac{1}{2}\sum\limits_{k=1}^{\left[\frac{m+1}{2}\right]-1}\sum\limits_{\alpha=0}^{2k-1}
\frac{(h\beta)^{2k-1-\alpha}(-1)^{\alpha}}{(2k-1-\alpha)!(\alpha+1)!}+
\frac{1}{2}\sum\limits_{k=\left[\frac{m+1}{2}\right]}^{m-1}
\sum\limits_{\alpha=0}^{m-2}
\frac{(h\beta)^{2k-1-\alpha}(-1)^{\alpha}}{(2k-1-\alpha)!(\alpha+1)!}+
$$
$$
+\frac{1}{2}\sum\limits_{k=\left[\frac{m+1}{2}\right]}^{m-1}
\sum\limits_{\alpha=m-1}^{m-2}
\frac{(h\beta)^{2k-1-\alpha}(-1)^{\alpha}}{(2k-1-\alpha)!(\alpha+1)!}.\eqno (4.39)
$$
Since the last equation is the identity with respect to $(h\beta)$
then from (4.39) equating the coefficients of the term
$(h\beta)^{2m-2}$ we obtain that for $m\geq 2$
$$
T=D_m(h\beta)*f_m(h\beta)=h.\eqno (4.40)
$$
Note that equality (4.40) for some $m$ was proved by directly
calculation of the convolution $D_m(h\beta)*f_m(h\beta)$ in [23].

Taking into account (4.40) for optimal coefficients (4.25) the following formula holds
$$
C_{\beta}=h+\sum\limits_{k=1}^{m-1}\left(a_k\lambda_k^{\beta}+b_k\lambda_k^{N-\beta}\right),\
\ \beta=\overline{1,N-1}.\eqno (4.41)
$$

Keeping in mind (4.40) after some simplifications from (4.39) we get
$$
\frac{e^{h\beta}}{2}\left[C_0+\frac{h}{e^h-1}+\sum\limits_{k=1}^{m-1}
\left(a_k\frac{\lambda_k}{e^h-\lambda_k}+b_k\frac{\lambda_k^N}{\lambda_ke^h-1}\right)-
\frac{e-1}{2e}\right]+
$$
$$
+\frac{e^{-h\beta}}{2}\left[-C_0-\frac{he^h}{1-e^h}+\sum\limits_{k=1}^{m-1}
\left(a_k\frac{\lambda_ke^h}{1-e^h\lambda_k}+b_k\frac{\lambda_k^Ne^h}{\lambda_k-e^h}\right)+
\frac{1}{2}\sum\limits_{\gamma=0}^NC_{\gamma}e^{h\gamma}+2d\right]+
$$
$$
+\frac{h(1+e^h)}{2(1-e^h)}-C_0\sum\limits_{k=1}^{m-1}\frac{(h\beta)^{2k-1}}{(2k-1)!}
-\sum\limits_{\ell=1}^{m-1}{h^{2\ell}}\sum\limits_{j=1}^{2\ell-1}
\frac{B_{2\ell-j}}{j!(2\ell-j)!}\beta^j+
$$
$$
+\sum\limits_{\ell=1}^{m-1}\frac{h^{2\ell-1}}{(2\ell-1)!}
\sum\limits_{j=0}^{2\ell-1}C_{2\ell-1}^j\beta^j
\sum\limits_{k=1}^{m-1}
a_k\frac{\lambda_k}{\lambda_k-1}\sum\limits_{i=0}^{2\ell-1}\frac{\Delta^i0^{2\ell-1-j}}
{(\lambda_k-1)^i}+
$$
$$
+\sum\limits_{\ell=1}^{m-1}\frac{h^{2\ell-1}}{(2\ell-1)!}
\sum\limits_{j=0}^{2\ell-1}C_{2\ell-1}^j\beta^j
\sum\limits_{k=1}^{m-1}
b_k\frac{\lambda_k^N}{1-\lambda_k}\sum\limits_{i=0}^{2\ell-1}
\left(\frac{\lambda_k}{1-\lambda_k}\right)^i\Delta^i0^{2\ell-1-j}+
$$
$$
+\frac{1}{2}\sum\limits_{k=[\frac{m+1}{2}]}^{m-1}\sum\limits_{\alpha=m-1}^{2k-1}
\frac{(h\beta)^{2k-1-\alpha}(-1)^{\alpha}}{(2k-1-\alpha)!\alpha!}
\sum\limits_{\gamma=0}^NC_{\gamma}(h\gamma)^{\alpha}+P_{m-2}(h\beta)=
$$
$$
=\frac{e^{h\beta}(1+e)}{4e}+\frac{e^{-h\beta}(1+e)}{4}-1+
\frac{1}{2}\sum\limits_{k=[\frac{m+1}{2}]}^{m-1}\sum\limits_{\alpha=m-1}^{2k-1}
\frac{(h\beta)^{2k-1-\alpha}(-1)^{\alpha}}{(2k-1-\alpha)!(\alpha+1)!}.
\eqno (4.42)
$$
As said above (4.42) is the identity with respect to $(h\beta)$.
From (4.42) equating the corresponding coefficients of
$e^{h\beta}$, $e^{-h\beta}$, $(h\beta)^{\alpha}$,
$\alpha=\overline{0,2m-3}$ one can get the system of linear
equations with respect to unknowns $a_k$, $b_k$
($k=\overline{1,m-1}$), $P_{m-2}(h\beta)$ and $d$.

Here we obtain the linear system for unknowns $a_k$, $b_k$
($k=\overline{1,m-1}$).

Equating the coefficients of $e^{h\beta}$ in the both
sides of (4.41) we get the following equation
$$
\sum\limits_{k=1}^{m-1}\left(a_k\frac{\lambda_k-\lambda_k^{N+1}}
{(e^h-\lambda_k)(1-\lambda_k)}+b_k\frac{\lambda_k-\lambda_k^{N+1}}
{(\lambda_ke^h-1)(1-\lambda_k)}\right)=0.\eqno (4.43)
$$
Further in (4.42) we consider the terms which consists of
$(h\beta)^{\alpha}$, $\alpha=\overline{m-1,2m-3}$ and we get the
equation
$$
\sum\limits_{\ell=[\frac{m+1}{2}]}^{m-1}\Bigg[-C_0\frac{(h\beta)^{2\ell-1}}{(2\ell-1)!}-
\sum\limits_{j=m-1}^{2\ell-1}(h\beta)^j\frac{h^{2\ell-j}B_{2\ell-j}}{j!\
(2\ell-j)!}+
\sum\limits_{j=m-1}
^{2\ell-1}\frac{(h\beta)^jh^{2\ell-j-1}}{(2\ell-1-j)!j!}\sum\limits_{k=1}^{m-1}
a_k\sum\limits_{i=0}^{2\ell-1}\frac{\lambda_k\Delta^i0^{2\ell-1-j}}
{(\lambda_k-1)^{i+1}}+
$$
$$+
\sum\limits_{j=m-1}
^{2\ell-1}(h\beta)^j\frac{h^{2\ell-j-1}}{(2\ell-1-j)!j!}\sum\limits_{k=1}^{m-1}
b_k\sum\limits_{i=0}^{2\ell-1}\frac{\lambda_k^{N+i}\Delta^i0^{2\ell-1-j}}
{(1-\lambda_k)^{i+1}}\Bigg]=0.\eqno (4.44)
$$

Now from equations (4.3) when $\alpha=0$ and (4.4) using
identities (4.26), (4.27), (4.28) taking into account (4.41) after
some simplifications for the coefficients $C_0$ and $C_N$ we
get the following expressions which are asserted in the theorem:\\
$$
\begin{array}{ll}
C_0=&\displaystyle\frac{e^h-1-h}{e^h-1}+\sum\limits_{k=1}^{m-1}
\Bigg(a_k\frac{\lambda_k(e^h-e)+\lambda_k^2(e-1)+\lambda_k^{N+1}(1-e^h)}
{(e-1)(1-\lambda_k)(e^h-\lambda_k)}+\\
&\displaystyle
+b_k\frac{\lambda_k^{N+1}(e^h-e)+\lambda_k^N(e-1)+\lambda_k(1-e^h)}
{(e-1)(\lambda_k-1)(\lambda_ke^h-1)}\Bigg),\\
\end{array}
\eqno (4.45)
$$
$$
\begin{array}{ll}
C_N=&\displaystyle\frac{e^hh+1-e^h}{e^h-1}+\sum\limits_{k=1}^{m-1}
\Bigg(a_k\frac{\lambda_k(e-e^{h+1})+\lambda_k^N(e^{h+1}-e^h)+\lambda_k^{N+1}(e^h-e)}
{(e-1)(1-\lambda_k)(e^h-\lambda_k)}+ \\
&
\displaystyle+b_k\frac{\lambda_k^{N+1}(e-e^{h+1})+\lambda_k^2(e^{h+1}-e^h)+\lambda_k(e^h-e)}
{(e-1)(1-\lambda_k)(1-\lambda_ke^h)}\Bigg).\\
\end{array}
\eqno (4.46)
$$
From (4.44), using (4.45),  grouping the coefficients of same
degrees of $(h\beta)$ and equating to zero that coefficients for
unknowns $a_k$ and $b_k$ we obtain the following $m-1$ linear
equations
$$
\begin{array}{l}
\sum\limits_{k=1}^{m-1}a_k\Bigg[\sum\limits_{l=1}^j\frac{h^{2l-2}}{(2l-2)!}\sum\limits_{i=0}^{2l-2}
\frac{\lambda_k\Delta^i0^{2l-2}}{(\lambda_k-1)^{i+1}}-
\frac{\lambda_k(e^h-e)+\lambda_k^2(e-1)+\lambda_k^{N+1}(1-e^h)}{(e-1)(\lambda_k-1)(\lambda_k-e^h)}\Bigg]+\\
+\sum\limits_{k=1}^{m-1}b_k\Bigg[\sum\limits_{l=1}^j\frac{h^{2l-2}}{(2l-2)!}\sum\limits_{i=0}^{2l-2}
\frac{\lambda_k^{N+i}\Delta^i0^{2l-2}}{(1-\lambda_k)^{i+1}}-
\frac{\lambda_k(1-e^h)+\lambda_k^N(e-1)+\lambda_k^{N+1}(e^h-e)}{(e-1)(\lambda_k-1)(\lambda_ke^h-1)}\Bigg]=\frac{e^h-1-h}{e^h-1}-\frac{h}{2},
j=\overline{1,[\frac{m}{2}]},
\end{array}
\eqno (4.47)
$$
$$
\begin{array}{l}
\sum\limits_{k=1}^{m-1}a_k\Bigg[\sum\limits_{l=1}^j\frac{h^{2l-1}}{(2l-1)!}\sum\limits_{i=0}^{2l-1}
\frac{\lambda_k\Delta^i0^{2l-1}}{(\lambda_k-1)^{i+1}}\Bigg]
+\sum\limits_{k=1}^{m-1}b_k\Bigg[\sum\limits_{l=1}^j\frac{h^{2l-1}}{(2l-1)!}\sum\limits_{i=0}^{2l-1}
\frac{\lambda_k^{N+i}\Delta^i0^{2l-1}}{(1-\lambda_k)^{i+1}}\Bigg],
j=\overline{1,[\frac{m-1}{2}]},
\end{array}
\eqno (4.48)
$$
Further from (4.3) when $\alpha=1,...,m-2$ using identities
(4.26), (4.27), (4.28) and the expression (4.46) for unknowns
$a_k$ and $b_k$ we have $m-2$ linear equations
$$
\begin{array}{ll}
\sum\limits_{k=1}^{m-1}a_k\Bigg[h^j\sum\limits_{i=0}^j\frac{\lambda_k^i-\lambda_k^{N+i}}{(1-\lambda_k)^{i+1}}\Delta^i0^j-
\sum\limits_{l=0}^{j-1}h^lC_j^l\sum\limits_{i=0}^l\frac{\lambda_k^{N+i}\Delta^i0^l}{(1-\lambda_k)^{i+1}}+
\frac{\lambda_k(e-e^{h+1})+\lambda_k^N(e^{h+1}-e^h)+\lambda_k^{N+1}(e^h-e)}{(e-1)(\lambda_k-1)(\lambda_k-e^h)}
\Bigg]+\\
+\sum\limits_{k=1}^{m-1}b_k\Bigg[h^j\sum\limits_{i=0}^j\frac{\lambda_k^{N+1}-\lambda_k}{(\lambda_k-1)^{i+1}}\Delta^i0^j-
\sum\limits_{l=0}^{j-1}h^lC_j^l\sum\limits_{i=0}^l\frac{\lambda_k\Delta^i0^l}{(\lambda_k-1)^{i+1}}+
\frac{\lambda_k^{N+1}(e-e^{h+1})+\lambda_k(e^{h}-e)+\lambda_k^2(e^{h+1}-e^h)}{(e-1)(\lambda_k-1)(\lambda_ke^h-1)}
\Bigg]=\\
\ \ \  \
=-\sum\limits_{l=1}^j\frac{j!B_{j+1-l}}{l!(j+1-l)!}h^{j+1-l}-\frac{e^hh+1-e^h}{e^h-1},\
\ j=\overline{1,m-2}.
\end{array}
\eqno (4.49)
$$

After some simplifications system of equations (4.43), (4.47),
(4.48), (4.49) we get the system which given in the assertion of
the theorem.

Theorem 4.6 is proved.\hfill $\Box$\\[0.2cm]

From theorem 4.6 when $m=2$ we get theorem 4.5.

For $m=3$ and $m=4$ from theorem 4.6 we we have the following
results.

\textbf{Corollary 4.1.} {\it The coefficients of optimal
quadrature formulas of the form (1.1) with the error functional
(1.2) and with equal spaced nodes in the space $W_2^{(3,2)}(0,1)$
are expressed by formulas}\\
$ C_0= \frac{e^h-1-h}{e^h-1}+\sum\limits_{k=1}^{2}
\Bigg(a_k\frac{\lambda_k(e^h-e)+\lambda_k^2(e-1)+\lambda_k^{N+1}(1-e^h)}
{(e-1)(1-\lambda_k)(e^h-\lambda_k)}+b_k\frac{\lambda_k^{N+1}(e^h-e)+\lambda_k^N(e-1)+\lambda_k(1-e^h)}
{(e-1)(\lambda_k-1)(\lambda_ke^h-1)}\Bigg),
$ \\
$ C_\beta=
h+\sum\limits_{k=1}^{2}\left(a_k\lambda_k^\beta+b_k\lambda_k^{N-\beta}\right),
\ \ \ \ \beta=\overline{1,N-1},\\
$\\
$ C_N= \frac{e^hh+1-e^h}{e^h-1}+\sum\limits_{k=1}^{2}
\Bigg(a_k\frac{\lambda_k(e-e^{h+1})+\lambda_k^N(e^{h+1}-e^h)+\lambda_k^{N+1}(e^h-e)}
{(e-1)(1-\lambda_k)(e^h-\lambda_k)}+
b_k\frac{\lambda_k^{N+1}(e-e^{h+1})+\lambda_k^2(e^{h+1}-e^h)+\lambda_k(e^h-e)}
{(e-1)(1-\lambda_k)(1-\lambda_ke^h)}\Bigg),
$ \\
\emph{where $a_k$ and $b_k$ ($k=\overline{1,2}$) are defined by
the following system of linear equations}
$$
\begin{array}{l}
\sum\limits_{k=1}^2a_k\frac{\lambda_k}{(\lambda_k-1)(\lambda_k-e^h)}+
\sum\limits_{k=1}^2b_k\frac{\lambda_k^{N+1}}{(\lambda_k-1)(\lambda_ke^h-1)}=\frac{h-2}{2(e^h-1)}+\frac{h}{(e^h-1)^2};\\
\sum\limits_{k=1}^2a_k\frac{\lambda_k^{N+1}}{(\lambda_k-1)(\lambda_k-e^h)}+
\sum\limits_{k=1}^2b_k\frac{\lambda_k}{(\lambda_k-1)(\lambda_ke^h-1)}=\frac{h-2}{2(e^h-1)}+\frac{h}{(e^h-1)^2};\\
\sum\limits_{k=1}^2a_k\frac{\lambda_k}{(\lambda_k-1)^2}+
\sum\limits_{k=1}^2b_k\frac{\lambda_k^{N+1}}{(\lambda_k-1)^2}=\frac{h}{12};\\
\sum\limits_{k=1}^2a_k\frac{\lambda_k^{N+1}}{(\lambda_k-1)^2}+
\sum\limits_{k=1}^2b_k\frac{\lambda_k}{(\lambda_k-1)^2}=\frac{h}{12};\\
\end{array}
$$
{\it here $\lambda_k$, $k=1,2$ are the roots of the polynomial
$$
\mathcal{P}_4(\lambda)=(1-e^{2h})(1-\lambda)^4-2(\lambda(e^{2h}+1)-e^h(\lambda^2+1))(h(1-\lambda)^2+\frac{h^3}{6}(1+4\lambda+\lambda^2))
$$
which $|\lambda_k|<1$.}

\textbf{Corollary 4.2.} {\it The coefficients of optimal
quadrature formulas of the form (1.1) with the error functional
(1.2) and with equal spaced nodes in the space $W_2^{(4,3)}(0,1)$
are expressed by formulas}\\
$ C_0= \frac{e^h-1-h}{e^h-1}+\sum\limits_{k=1}^{3}
\Bigg(a_k\frac{\lambda_k(e^h-e)+\lambda_k^2(e-1)+\lambda_k^{N+1}(1-e^h)}
{(e-1)(1-\lambda_k)(e^h-\lambda_k)}+b_k\frac{\lambda_k^{N+1}(e^h-e)+\lambda_k^N(e-1)+\lambda_k(1-e^h)}
{(e-1)(\lambda_k-1)(\lambda_ke^h-1)}\Bigg),
$ \\
$ C_\beta=
h+\sum\limits_{k=1}^{3}\left(a_k\lambda_k^\beta+b_k\lambda_k^{N-\beta}\right),
\ \ \ \ \beta=\overline{1,N-1},\\
$\\
$ C_N= \frac{e^hh+1-e^h}{e^h-1}+\sum\limits_{k=1}^{3}
\Bigg(a_k\frac{\lambda_k(e-e^{h+1})+\lambda_k^N(e^{h+1}-e^h)+\lambda_k^{N+1}(e^h-e)}
{(e-1)(1-\lambda_k)(e^h-\lambda_k)}+
b_k\frac{\lambda_k^{N+1}(e-e^{h+1})+\lambda_k^2(e^{h+1}-e^h)+\lambda_k(e^h-e)}
{(e-1)(1-\lambda_k)(1-\lambda_ke^h)}\Bigg),
$ \\
\emph{where $a_k$ and $b_k$ ($k=\overline{1,3}$) are defined by
the following system of linear equations}
$$
\begin{array}{l}
\sum\limits_{k=1}^3a_k\frac{\lambda_k}{(\lambda_k-1)(\lambda_k-e^h)}+
\sum\limits_{k=1}^3b_k\frac{\lambda_k^{N+1}}{(\lambda_k-1)(\lambda_ke^h-1)}=\frac{h-2}{2(e^h-1)}+\frac{h}{(e^h-1)^2};\\
\sum\limits_{k=1}^3a_k\frac{\lambda_k^{N+1}}{(\lambda_k-1)(\lambda_k-e^h)}+
\sum\limits_{k=1}^3b_k\frac{\lambda_k}{(\lambda_k-1)(\lambda_ke^h-1)}=\frac{h-2}{2(e^h-1)}+\frac{h}{(e^h-1)^2};\\
\sum\limits_{k=1}^3a_k\frac{\lambda_k}{(\lambda_k-1)^2}+
\sum\limits_{k=1}^3b_k\frac{\lambda_k^{N+1}}{(\lambda_k-1)^2}=\frac{h}{12};\\
\sum\limits_{k=1}^3a_k\frac{\lambda_k^{N+1}}{(\lambda_k-1)^2}+
\sum\limits_{k=1}^3b_k\frac{\lambda_k}{(\lambda_k-1)^2}=\frac{h}{12};\\
\sum\limits_{k=1}^3a_k\frac{\lambda_k}{(\lambda_k-1)^3}+
\sum\limits_{k=1}^3b_k\frac{\lambda_k^{N+2}}{(1-\lambda_k)^3}=-\frac{h}{24};\\
\sum\limits_{k=1}^3a_k\frac{\lambda_k^2-\lambda_k^{N+2}}{(1-\lambda_k)^3}+
\sum\limits_{k=1}^3b_k\frac{\lambda_k^{N+1}-\lambda_k}{(\lambda_k-1)^3}=0;\\
\end{array}
$$
{\it here $\lambda_k$, $k=1,2,3$ are the roots of the polynomial
$$
\begin{array}{ll}
\mathcal{P}_6(\lambda)=&(1-e^{2h})(1-\lambda)^6-2(\lambda(e^{2h}+1)-
         e^h(\lambda^2+1))\times\\
 & \times\bigg(h(1-\lambda)^4+\frac{h^3}{6}(1-\lambda)^2(1+4\lambda+\lambda^2)+
 \frac{h^5}{120}(1+26\lambda+66\lambda^2+26\lambda^3+\lambda^4)\bigg)
\end{array}
$$
which $|\lambda_k|<1$.}

\section{Numerical results}

In this section we give the numerical results which confirm the
theoretical results obtained in the section 4.

Taking into account (2.9) from (2.21) for the norm of the error
functional of optimal quadrature formulas in the space
$W_2^{(m,m-1)}(0,1)$ we obtain
$$
\begin{array}{ll}
\left\| \ell(x)\right\|^2 = &\displaystyle
(-1)^m\Bigg[\sum\limits_{\beta = 0}^N \sum\limits_{\gamma  = 0}^N
C_\beta  C_\gamma  \,
\frac{\mathrm{sing}(h\beta-h\gamma)}{2}\left(\frac{e^{h\beta-h\gamma}-e^{h\gamma-h\beta}}{2}-\sum\limits_{k=1}^{m-1}
\frac{(h\beta-h\gamma)^{2k-1}}{(2k-1)!}\right) - \\
 &\displaystyle  - 2\sum\limits_{\beta  = 0}^N C_\beta
\left(\frac{e^{h\beta}+e^{h\beta}+e^{1-h\beta}+e^{h\beta-1}-4}{4}-\sum\limits_{k=1}^{m-1}\frac{(h\beta)^{2k}+(1-h\beta)^{2k}}{2\cdot
(2k)!}\right)+ \\
&\displaystyle+\frac{e^2-2e-1}{2e}-\sum\limits_{k=1}^{m-1}\frac{1}{(2k+1)!}\Bigg],
\end{array}
 \eqno   (5.1)
$$
In computations of optimal coefficients we need the roots
$\lambda_k$, $k=\overline{1,m-1}$, $|\lambda_k|<1$ of the
polynomial defined by equality (4.12)
$$
\mathcal{P}_{2m-2}(\lambda)=\sum\limits_{s=0}^{2m-2}p_s^{(2m-2)}\lambda^s=(1-e^{2h})(1-\lambda)^{2m-2}-
2(\lambda(e^{2h}+1)-e^h(\lambda^2+1))\times
$$
$$
\times
\left[h(1-\lambda)^{2m-4}+\frac{h^3(1-\lambda)^{2m-6}}{3!}E_2(\lambda)+...+
\frac{h^{2m-3}E_{2m-4}(\lambda)}{(2m-3)!}\right],\eqno (5.2)
$$
where $E_{2k}(\lambda)=\sum\limits_{s=0}^{2k}e_s\lambda^s$ is the
Euler-Frobenius polynomial of degree $2k$ and for the coefficients
of this polynomial the formula $e_s=\sum\limits_{j=0}^s(-1)^j{2k+2
\choose j} (s+1-j)^{2k+1}$ is used which was given by Euler.

For convenience the absolute value of the difference (1.4) of the
quadrature formula (1.1) we denote by $|R(\varphi)|$. Then by
Cauchy-Scwartz inequality we have
$$
|R(\varphi)|\leq \|\varphi(x)|W_2^{(m,m-1)}(0,1)\|\cdot
\|\ell(x)|W_2^{(m,m-1)*}(0,1)\|.\eqno (5.3)
$$
We consider the cases $m=1,\ 2,\ 3,\ 4$ and $N=10,\ 50,\ 100$.

 \textbf{The case $m=1$.} In the space $W_2^{(1,0)}(0,1)$ using theorem 4.4 and (5.1), (5.3) for the error of
 optimal quadrature formula (1.1) we have
 $$
\begin{array}{ll}
N=10:\  &|R(\varphi)|\leq \|\varphi(x)|W_2^{(1,0)}(0,1)\|\cdot 0.02886;\\
N=50:\  &|R(\varphi)|\leq \|\varphi(x)|W_2^{(1,0)}(0,1)\|\cdot 0.00577;\\
N=100:\ &|R(\varphi)|\leq \|\varphi(x)|W_2^{(1,0)}(0,1)\|\cdot 0.00289.\\
\end{array}
 $$

 \textbf{The case $m=2$.} In the space $W_2^{(2,1)}(0,1)$ using theorem 4.5 and (5.1), (5.2), (5.3) for the error of
 optimal quadrature formula (1.1) we have
 $$
\begin{array}{ll}
N=10:\  &|R(\varphi)|\leq \|\varphi(x)|W_2^{(2,1)}(0,1)\|\cdot 0.000424;\\
N=50:\  &|R(\varphi)|\leq \|\varphi(x)|W_2^{(2,1)}(0,1)\|\cdot 0.00001534;\\
N=100:\ &|R(\varphi)|\leq \|\varphi(x)|W_2^{(2,1)}(0,1)\|\cdot 0.37802\times 10^{-5}.\\
\end{array}
 $$

\textbf{The case $m=3$.} In the space $W_2^{(3,2)}(0,1)$ using
corollary 4.1 and (5.1), (5.2), (5.3) for the error of optimal
quadrature formula (1.1) we have
 $$
\begin{array}{ll}
N=10:\  &|R(\varphi)|\leq \|\varphi(x)|W_2^{(3,2)}(0,1)\|\cdot 0.0000108;\\
N=50:\  &|R(\varphi)|\leq \|\varphi(x)|W_2^{(3,2)}(0,1)\|\cdot 0.5643\times 10^{-7};\\
N=100:\ &|R(\varphi)|\leq \|\varphi(x)|W_2^{(3,2)}(0,1)\|\cdot 0.6435\times 10^{-8}.\\
\end{array}
 $$

\textbf{The case $m=4$.} In the space $W_2^{(4,3)}(0,1)$ using
corollary 4.2 and (5.1), (5.2), (5.3) for the error of optimal
quadrature formula (1.1) we have
 $$
\begin{array}{ll}
N=10:\  &|R(\varphi)|\leq \|\varphi(x)|W_2^{(4,3)}(0,1)\|\cdot 0.5051\times 10^{-6};\\
N=50:\  &|R(\varphi)|\leq \|\varphi(x)|W_2^{(4,3)}(0,1)\|\cdot 0.3854\times 10^{-9};\\
N=100:\ &|R(\varphi)|\leq \|\varphi(x)|W_2^{(4,3)}(0,1)\|\cdot 0.1821\times 10^{-10}.\\
\end{array}
 $$

\section*{ Acknowledgements}

The authors are very thankful to professor Erich Novak for
discussion of the results of this paper.

The final part of this work was done in the Friedrich-Schiller
University of Jena, Germany. The second author thanks the DAAD for
scholarship. Furthermore the second author also thanks
IMU/CDE-program for travel support to the Friedrich-Schiller
University of Jena, Germany.

\end{document}